\documentclass[10pt]{amsart}
\usepackage{amsfonts}
\usepackage{color}
\usepackage[square,compress,comma, numbers]{natbib}
\usepackage[colorlinks=true, citecolor=blue, linkcolor=blue]{hyperref}
\allowdisplaybreaks[4]
\usepackage{amssymb}
\usepackage{amsmath}

\definecolor{c20}{rgb}{0.,0.7,0.}
\definecolor{c30}{rgb}{0.,0.,1.}
\definecolor{c40}{rgb}{1,0.1,0.7}
\definecolor{c50}{rgb}{1,0,0}
\definecolor{c60}{rgb}{1,0.9,0.1}

\allowdisplaybreaks[4]

\def\cLa#1{\textcolor{c50}{#1}}
\def\cLa#1{#1}

\def\kk#1{\textcolor{black}{#1}}

\newcommand{\kb}[1]{\boldsymbol{#1}}
\newcommand{\vk}[1]{\kb{#1}}

\newcommand{\abs}[1]{\left\lvert #1 \right\rvert}

\newcommand{\E}[1]{\mathbb{E}\left\{#1\right\}}

\newcommand{\pk}[1]{\mathbb{P} \left \{#1 \right \} }

\newcommand{\R}{\mathbb{R}}

\newcommand{\N}{\mathbb{N}}

\newcommand{\ldot}{,\ldots,}

\newcommand{\BQN}{\begin{eqnarray}}
\newcommand{\EQN}{\end{eqnarray}}
\newcommand{\BQNY}{\begin{eqnarray*}}
\newcommand{\EQNY}{\end{eqnarray*}}

\newcommand{\BS}{\begin{sat}}
\newcommand{\ES}{\end{sat}}
\newcommand{\BT}{\begin{theo}}
\newcommand{\ET}{\end{theo}}
\newcommand{\BL}{\begin{lem}}
\newcommand{\EL}{\end{lem}}
\newcommand{\BK}{\begin{korr}}
\newcommand{\EK}{\end{korr}}

\newcommand{\BD}{\begin{de}}
\newcommand{\ED}{\end{de}}

\newcommand{\BIT}{\begin{itemize}}
\newcommand{\EIT}{\end{itemize}}
\newcommand{\BDI}{\begin{description}}
\newcommand{\EDI}{\end{description}}

\newcommand{\BRM}{\begin{remarks}}
\newcommand{\ERM}{\end{remarks}}

\newcommand{\BEL}{\begin{lem}}
\newcommand{\EEL}{\end{lem}}

\newtheorem{theo}{Theorem}[section]
\newtheorem{sat}[theo]{Proposition}
\newtheorem{de}[theo]{Definition}
\newtheorem{lem}[theo]{Lemma}

\newtheorem{korr}[theo]{Corollary}
\newtheorem{remark}[theo]{Remark}
\newtheorem{remarks}[theo]{Remarks}
\newtheorem{prop}[theo]{Proposition}

\newcommand{\nelem}[1]{{Lemma \ref{#1}}}

\newcommand{\netheo}[1]{{Theorem \ref{#1}}}

\newcommand{\proofprop}[1]{\textsc{\bf Proof of Proposition} \ref{#1}:}

\newcommand{\COM}[1]{}

\newcommand{\QED}{\hfill $\Box$}

\def\rw{\rightarrow}

\def\IF{\infty}

\def\LT{\left}
\def\RT{\right}

\def\rw{\rightarrow}

\def\vn{\varepsilon}
\def\Var{\text{Var}}
\def\aalpha{\vk{\alpha}}
\def\bbeta{\vk{\beta}}
\def\ta{\vk{a}}

\def\x{\vk{x}}
\def\y{\vk{y}}
\def\z{\vk{z}}
\def\t{\vk{t}}

\def\s{\vk{s}}

\def\l{\vk{l}}
\def\Sn{\mathcal{S}_{n}}

\def\MNB{\mathcal{N}^c}
\newcommand{\pu}[1]{\mathbf{P}_u\left(#1\right)}

\def\l{\vk{l}}
\def\j{\vk{j}}

\def\Du{\mathcal{D}_{u}}

\def\tt{\tilde{\t}}
\def\bt{\bar{\t}}
\def\tz{\tilde{\z}}
\def\ll{{\langle}}
\def\r{\rangle}

\def\x{{\bf{x}}}
\def\y{{\bf{y}}}
\def\o{{\bf{o}}}
\def\laweq{\buildrel d \over =}
\begin{document}

\title[ Extremes of Gaussian random fields with non-additive dependence structure  ]{Extremes of Gaussian random fields  with non-additive dependence structure }

\author{Long Bai}
\address{Long Bai,
Department of Statistics and Actuarial Science , Xi'an Jiaotong-Liverpool University,
Suzhou 215123
China}
\email{long.bai@xjtlu.edu.cn}

\author{Krzysztof D\c{e}bicki}

\address{Krzysztof D\c{e}bicki,
	Mathematical Institute,
	University of Wroclaw, pl. Grunwaldzki 2/4, 50-384 Wroclaw, Poland}
\email{Krzysztof.Debicki@math.uni.wroc.pl}

\author{Peng Liu}
\address{Peng Liu, Department of Mathematical Sciences,\\
	Wivenhoe Park, Colchester CO4 3SQ, UK
}
\email{peng.liu@essex.ac.uk}

\bigskip

\date{\today}
 \maketitle

{\bf Abstract:}
We derive exact asymptotics
of
\[
\pk{\sup_{\t\in {\mathcal{A}}}X(\t)>u},\ {\rm as}\ u\to\infty,
\]
for
a centered Gaussian field $X({\t}),\ {\t}\in \mathcal{A}\subset\R^n$, $n>1$
with continuous sample paths a.s.,
for which
$\arg \max_{\t\in {\mathcal{A}}} Var(X(\t))$
is a Jordan set with
finite and positive Lebesque measure of dimension $k\le n$ and
\kk{its dependence structure
is not necessarily locally stationary.}
Our findings are applied to deriving the asymptotics of tail probabilities
related to  performance tables and chi processes where the covariance structure is
not locally stationary.



\section{Introduction}

Let $X({\t}),\ {\t}\in \R^n$, $n>1$ be a centered Gaussian field with continuous sample paths a.s..
Due to its importance in the extreme value theory of stochastic processes, statistics
and applied probability, the distributional properties of
\begin{eqnarray}
\sup_{\t\in {\mathcal{A}}}X(\t),\label{sup}
\end{eqnarray}
with a bounded set $\mathcal{A}\subset \R^n$, were vastly investigated.
While the exact distribution of (\ref{sup})
is known only for very particular processes, the asymptotics of
\begin{eqnarray}
\pk{\sup_{\t\in {\mathcal{A}}}X(\t)>u}\label{tail},
\end{eqnarray}
as $u\to\infty$ was intensively analyzed;
see the seminal monographs \cite{Pit96, AdlerTaylor, Lif13}.
As advocated  therein, the set
of points that maximize the variance $\mathcal{M}^\star:=\arg \max_{\t\in {\mathcal{A}}} Var(X(\t))$
plays an important role in the form of the exact asymptotics of (\ref{tail}).
The best understood cases cover the situation when
(i) $v_n(\mathcal{M}^\star)\in (0,\infty)$,
where $v_n$ denotes the Lebesgue measure on $\R^n$
and the field $X(\t)$ is homogeneous 
on $\mathcal{M}^\star$
or (ii) the set $\mathcal{M}^\star$
consists of separate points.
In case (i) one can argue that
\[\pk{\sup_{\t\in {\mathcal{A}}}X(\t)>u}\sim \pk{\sup_{\t\in {\mathcal{M}^\star}}X(\t)>u}\ {\rm as}\ u\to\infty.\]
For the intuitive description of case (ii), suppose that $\mathcal{M}^\star=\{\t^\star\}$ and $Var(X(\t^\star))=1$.
Then, the play between local behaviour of standard deviation and correlation function
in the neighbourhood of $\mathcal{M}^\star$ influences the asymptotics, which
takes the form
\begin{eqnarray}
\pk{\sup_{\t\in {\mathcal{A}}}X(\t)>u}\sim f(u)\pk{X(\t^\star)>u}\ {\rm as}\ u\to\infty.\label{single}
\end{eqnarray}
An applicable 
assumption
under which one can get the exact asymptotics as given in (\ref{single}),
is that
in the neighbourhood of $\t^\star$, both the standard deviation and correlation function
of $X(\t)$ factorizes according to the following additive form
\begin{eqnarray}
1- \sigma(\t) \sim \sum_{j=1}^3 g_j( \bt^\star_j- \bt_j)\label{as.1},\ \
1-Corr(\s,\t)\sim \sum_{j=1}^3 h_j(\overline{\s}_j-\bt_j)
\end{eqnarray}
as $\s,\t \to \t^\star$, where
the coordinates of $\R^n$ are split onto disjoint sets $\Lambda_1,\Lambda_2,\Lambda_3$
with $\{1,...,n\}=\Lambda_1\cup\Lambda_2\cup\Lambda_3$,
$\bt_j=(t_{i})_{i\in \Lambda_j},$  $j=1,2,3$
for $\t \in \R^n$ and
$g_j,h_j$ are some homogeneous functions (see  (\ref{G}))
such that
\begin{eqnarray}
\lim_{\bt_1\to \overline{0}_1}\frac{g_1(\bt_1)}{ h_1(\bt_1)}=0,
\
\lim_{\bt_2\to \overline{0}_2}\frac{g_2(\bt_2)}{ h_2(\bt_2)}\in (0,\infty),
\
\lim_{\bt_3\to \overline{0}_3}\frac{g_3(\bt_3)}{ h_3(\bt_3)}=\infty.\label{lims}
\end{eqnarray}
Importantly, under (\ref{as.1})-(\ref{lims})
\[f(u)=f_1(u)f_2(u)f_3(u)\]
factorizes too and:
\\
$\diamond$
in the direction of coordinates $\Lambda_1$ the standard deviation function is relatively flat with comparison
to the correlation function. Then,
with respect of coordinates  $\Lambda_1$,
some substantial neighbourhood of $\mathcal{M}$
contributes to the asymptotics
and $f_1(u)\to\infty$ as $u\to\infty$.\\
$\diamond$
in the direction of coordinates  $\Lambda_2$ the standard deviation function is comparable to the correlation function.
Then, with respect of coordinates  $\Lambda_2$, some relatively {\it small} neighbourhood of $\mathcal{M}$
is important for the asymptotics
and $f_2(u)\to \mathcal{P}\in (1,\infty)$ as $u\to\infty$.\\
$\diamond$ in the direction of coordinates  $\Lambda_3$
the standard deviation function decreases relatively fast
with comparison to the correlation. Then,
with respect of coordinates  $\Lambda_3$, only the sole optimizer $t^\star$ is
responsible for the asymptotics and $f_3(u)\to 1$ as $u\to\infty$.
We refer to Piterbarg \cite{Pit96}[Chapter 8] for the details.
\\
Much less is known on the mixed cases, when
set $\mathcal{M}^\star$ is a more general subset of $\mathcal{A}$
and/or
the local dependance structure of the analyzed process doesn't
factorize
according to the additive structure
as in (\ref{as.1})-(\ref{lims}).
The exemptions that are available in the literature  were analyzed separately
and cover some particular
cases as in \cite{PiP81,Liu16,Liu18,DHL16,Chan2006,Adler1986}.
We would like to point at a notable recent contribution by Piterbarg \cite{Pit21},
which 
deals with analysis
of high excursion probability for centered Gaussian fields on finite dimensional manifold
where $\mathcal{M}^\star$ is a smooth submanifold. In this intuitively presented
work, under the assumption that the correlation function of $X$
is locally homogeneous,
three scenarios for $\mathcal{M}^\star\varsubsetneq \mathcal{A}$ are worked out:
(1) {\it stationary like case}, (2) {\it transition case}
and (3) {\it Talagrand case};
which under the notation
(\ref{as.1})-(\ref{lims}) correspond to $\Lambda_2=\Lambda_3=\emptyset$ for (1),
$\Lambda_1=\Lambda_3=\emptyset$ for (2), $\Lambda_1=\Lambda_2=\emptyset$ for (3).

In view of the considered in our paper examples and transparency of the presentation of
the results, we  work on Euclidean space in this contribution.
We derive a unified result that allows to get exact asymptotics for the class
of centered Gaussian fields for which we allow that $\mathcal{M}^\star$ is a $k_0\le n$ dimensional
bounded Jordan set and the dependence structure of the entire field in the neighbourhood of $\mathcal{M}^\star$
doesn't necessarily decompose 
as in (\ref{as.1})-(\ref{lims}).
In comparison to \cite{Pit21},
we allow mixed scenarios where all sets $\Lambda_1,\Lambda_2,\Lambda_3$ can be  nonempty
at the same time.
Besides, we suppose that $X$ is locally stationary
only with respect to coordinates of stationary like direction (see assumption {\bf A1});
this relaxation is particularly important for the examples that are worked out in Section \ref{s.GUE} and Section \ref{s.chi}.

One of the motivations for this contribution is
the analysis of asymptotic properties of
\BQN\label{aim20}
\pk{D_n^{\alpha}>u}:=\pk{\sup_{\t\in\mathcal{S}_{n}}Z^{\alpha}(\t)>u}, \ {\rm as}\ u\to\infty,
\EQN
where
$\t=(t_1\ldot t_n)$, $\mathcal{S}_{n}=\{\t\in \R^n: 0\leq t_1\leq\cdots\leq t_n\leq 1\}$,
\BQNY
Z^{\alpha}(\t)=\sum_{i=1}^{n+1}a_i\LT(B^{\alpha}_{i}(t_i)-B^{\alpha}_{i}(t_{i-1})\RT)
\EQNY
with $t_0=0,t_{n+1}=1$
and $B^{\alpha}_{i},\ i=1\ldot n+1$ are independent fractional Brownian motions with Hurst index $\alpha/2\in(0,1)$.
This random variable plays an important role in many areas of probability theory.
In particular, for
${\alpha}_{i}\equiv 1$
it is strongly related with the notion of the performance table
and it also appears as a limit
in problems describing queues in series, totally asymmetric exclusion processes
or oriented percolation \cite{Bar01,Oco02, GlW91}.
If additionally $a_i\equiv 1$, then
$D_n^{1}$ has the same law as the largest eigenvalue of an n-dimensional
GUE (Gaussian Unitary Ensamble) matrix \cite{GTW01}.
However,
if $\alpha=1$ but $a_i$ are not all the same,
then the size of $\mathcal{M}^\star$ depends on the number of maximal $a_i$ and the correlation
structure of the entire field is not locally homogeneous.
Application of Theorem \ref{MThm2} in Section 2 allows to derive exact asymptotics of (\ref{aim20}) as $u\to\infty$
for $\alpha\in(0,2)$; see Proposition \ref{ThmM1}.

Another illustration of the applicability of Theorem \ref{MThm2}
deals with the extremes of the class of {\it chi} processes $\chi(t),t\ge0$, defined for
given $\vk{X}(t)=(X_1(t),\cdots, X_n(t)), t\ge 0$ where $X_i(t)$ for $i=1,...,n$ are mutually independent,
as
\[
\chi(t):=\sqrt{\sum_{i=1}^n X_i^2(t)}, \ t\ge0.
\]
Due to their importance in statistics, asymptotic properties of high excursions of chi processes attracted
substantial interest.
We refer to the classical work by Lindgren \cite{Lindgren1980a} 
and more recent contributions \cite{Pit21, Pitchi1994,EnkelejdJi2014Chi,PL2015,LJ2017,Bai2021}
that deal with non-stationary or noncentered cases.
In Section \ref{s.chi} we apply Theorem \ref{MThm2}
to the analysis of the asymptotics of tail distribution of high exceedances of $\chi(t)$
for a model, where the covariance structure of $X_i$ is not locally stationary; see Proposition \ref{ex02}.

 The structure of the rest of the paper is organized as follows. The proofs of Theorem \ref{MThm2},
 Proposition \ref{ThmM1} and Proposition \ref{ex02} are given in Sections 4-6
 respectively while the proofs of some auxiliary results are postponed to the Appendix.

\section{Main result}\label{s.m}

 Let $X(\t),\ \t\in \mathcal{A}$  be an $n$-dimensional centered Gaussian field with  continuous trajectories,
 variance function $\sigma^2(\t)$ and correlation function $r(\s,\t)$, where $\mathcal{A}$ is a \kk{bounded} set in $\R^n$.
Suppose that the maximum of variance function $\sigma^2(\t)$
over $\mathcal{A}$ is attained on a \kk{Jordan} subset of $\mathcal{A}$.
Without loss of generality, we assume that $\max_{\t\in \mathcal{A}} \sigma^2(\t)=1$ and we denote
$\mathcal{M}^*:=\{\t\in \mathcal{A}: \sigma^2(\t)=1 \}$.

Throughout this paper, all the operations on vectors are meant componentwise. For instance, for any given $\x = (x_1\ldot x_n)\in\R^n$ and $\y = (y_1\ldot y_n)\in \R^n$, we write $\x\y=(x_1y_1\ldot x_ny_n)$, $1/\x=(1/x_1,\dots, 1/x_n)$ for $x_i> 0, i=1,\dots, n$, and $\x^{\y}=(x_1^{y_1},\dots, x_n^{y_n})$ for $x_i, y_i\geq 0, i=1,\dots, n$. Moreover, we say that
$\x\geq \y$ for $x_i\geq y_i, i=1,\dots, n$.

Suppose that the coordinates of $\R^n$ can be exclusively
split onto four disjoint sets $\Lambda_i, i=0,1,2,3$
with
$k_i=\#\cup_{j=0}^i\Lambda_j, ~i=0,1,2,3$ (implying that
$1\leq k_0\leq k_1\leq k_2\leq k_3 $ with $k_3=n$)
and
\BQNY
 \tt:=(t_i)_{i\in \Lambda_0},~\bt_j:=(t_{i})_{i\in \Lambda_j}, \quad j=1,2,3
\EQNY in such a way that
$\mathcal{M}^*=\{\t\in\mathcal{A}: t_{i}=0, i\in \cup_{j=1,2,3}\Lambda_j\}$.
Note that $\mathcal{M}^*=\mathcal{A}$ if $\cup_{j=1,2,3}\Lambda_j=\emptyset$.
Sets $\Lambda_1,\Lambda_2,\Lambda_3$ play similar role to the described
in the Introduction (see {\bf A2} below), while $\Lambda_0$ relates to
$\mathcal{M}^\star$ by
$\mathcal{M}:=\{\tt:\t\in\mathcal{A}, t_{i}=0, i\in \cup_{j=1,2,3}\Lambda_j\}\subset\mathbb{R}^{k_0}$.

Suppose that $\mathcal{M}$ is Jordan measurable with $v_{k_0}(\mathcal{M})\in (0,\IF)$,
where
$v_{k_0}$ denotes the Lebesgue measure on $\mathbb{R}^{k_0}$,
and  $\{(t_1,\dots, t_n): \tt\in\mathcal{M}, t_i\in [0,\epsilon), i\in \cup_{j=1,2,3}\Lambda_j \}\subseteq \mathcal{A}\subseteq \{(t_1,\dots, t_n): \tt\in\mathcal{M}, t_i\in [0,\infty), i\in \cup_{j=1,2,3}\Lambda_j \}$
for some $\vn\in(0,1)$
small enough.
Further, we shall impose the following assumptions on the standard deviation and correlation functions of $X$:\\

{\bf A1}: There exists a centered Gaussian random  field $W(\t),\ \t\in[0,\IF)^n$ with continuous sample paths and a positive continuous vector-valued function  $\vk{a}(\tz)=(a_1(\tz)\ldot a_n(\tz)),\ \tz=(z_i)_{i\in\Lambda_0}\in \mathcal{M}$ satisfying
\BQN\label{az}
\inf_{i=1,\dots ,n}\inf_{\tz\in \mathcal{M}}a_i(\tz)>0
\EQN
 such that
\BQN\label{Rr}
\lim_{\delta\rw 0}\sup_{\z\in\mathcal{M}^*}\underset{\abs{\s-\z},\abs{\t-\z}\leq \delta}{\sup_{\s,\t\in \mathcal{A}}}\abs{\frac{1-r(\s,\t)}{\E{\LT(W(\vk{a}(\tz)\s)-W(\vk{a}(\tz)\t)\RT)^2}}-1}=0,
\EQN
where the increments of $W$ {are} homogeneous if we fix both $\bt_2$ and $ \bt_3$, and there exists a vector $\aalpha=(\alpha_1,\dots, \alpha_n)$ with $\alpha_i\in (0,2],1\leq i\leq n$ such that for any $u>0$
\BQN\label{W}
\E{\LT(W(u^{-2/\vk{\alpha}}\s)-W(u^{-2/\vk{\alpha}}\t)\RT)^2}=u^{-2}\E{\LT(W(\s)-W(\t)\RT)^2}.
\EQN
Moreover, there exist $d>0, \mathcal{Q}_i>0,i=1,2$ such that  for any $\s,\t\in \mathcal{A}$ and $|\s-\t|<d$
\BQN\label{boundrr1}
\mathcal{Q}_1\sum_{i\in \cup_{j=0,1}\Lambda_j} \abs{s_i-t_i}^{\alpha_i}\leq 1-r(\s,\t) \leq \mathcal{Q}_2\sum_{i=1}^n \abs{s_i-t_i}^{\alpha_i}
\EQN
Further, suppose that for $\s,\t\in\mathcal{A}$ and $\s\neq \t$
\BQN\label{reql}
r(\s,\t)<1.
\EQN

{\bf A2}:  Assume that
\BQN\label{Var}
\lim_{\delta\rw 0}\sup_{\z\in\mathcal{M}^*}\underset{\abs{\z-\t}\leq\delta}{\sup_{\t \in \mathcal{A}}}
\abs{\frac{ 1- \sigma(\t)}{\sum_{j=1}^3p_j(\tz)g_j(\bt_j)}- 1}=0,
\EQN
where $p_j(\tt),\ \tt\in[0,\IF)^{k_0}, j=1,2,3$ are  positive continuous functions and
$g_j(\bt_j),\bt_j\in\mathbb{R}^{k_j-k_{j-1}}, j=1,2,3$, are  {continuous} functions satisfying
$g_i(\bt_i)>0, \bt_j\neq \overline{\bf 0}_j, j=1,2,3.$ Moreover, we shall assume the following homogeneous property on $g_j$'s:
for any $u>0$ and some $\bbeta_j=(\beta_{i})_{i\in \Lambda_j},~ j=1,2,3$ with $\beta_i>0, i\in \cup_{j=1,2,3}\Lambda_j, $
\BQN\label{G}
u g_j(\bt_j)&=&g_j(u^{1/{ \bbeta_{j}}}\bt_{j}), \quad j=1,2,3.
\EQN
Moreover, with $\aalpha_j=(\alpha_{i})_{i\in \Lambda_j}, \quad j=1,2,3$,
\begin{align}
\aalpha_1<\bbeta_1, \aalpha_2=\bbeta_2~ \text{and}~ \aalpha_3>\bbeta_3.
\end{align}

We next display the main result of this paper.
To the end of this paper $\Psi(\cdot)$ denotes tail distribution of the standard normal random variable.
\BT\label{MThm2}
Suppose that  $X(\t),\ \t\in \mathcal{A}$ is a  $n$-dimensional centered Gaussian random field satisfying  {\bf A1-A2}.
Then, as $u\to\infty$,
\BQNY
\pk{\sup_{\t\in\mathcal{A}}X(\t)>u}\sim
C u^{\sum_{i\in \Lambda_0\cup\Lambda_1}\frac{2}{\alpha_i}-\sum_{i\in \Lambda_1}\frac{2}{\beta_i}}\Psi(u),
\EQNY
where
\BQNY
C=  \int_{\mathcal{M}}\left(\mathcal{H}_W^{p_2(\tz)g_2(\vk{a}_2^{-1}(\tz)\bt_2)} \LT(\prod_{i\in \Lambda_0\cup\Lambda_1}a_i(\tz)\RT) \int_{\bt_1\in[0,\IF)^{k_1-k_0}} e^{-p_1(\tz)g_1(\bt_1)}d\bt_1\right) d\tz \in (0,\IF),
\EQNY
with $\vk{a}_2(\tz)=(a_i(\tz))_{i\in \Lambda_2}$ and
\BQNY
\mathcal{H}_W^{p_2(\tz)g_2(\vk{a}_2^{-1}(\tz)\bt_2)}= \lim_{\lambda\rw\IF}\frac{1}{\lambda^{k_1}}\E{\sup_{t_i\in[0,\lambda], i\in \Lambda_0\cup\Lambda_1\cup\Lambda_2; t_i=0, i\in \Lambda_3}e^{ \sqrt{2}W(\t)-\sigma^2_W(\t)-p_2(\tz)g_2(\vk{a}_2^{-1}(\tz)\bt_2)}}.
\EQNY
\ET
\begin{remark}\label{remark1}
The result in Theorem \ref{MThm2} is also valid if some of $\Lambda_i, i=0,1,2,3$ are empty sets.
\end{remark}
 Next, let us consider  a special case of Theorem \ref{MThm2}. Suppose that
 \begin{align}\label{ai}
 a_i(\tz)\equiv a_i, \tz\in \mathcal{M},~i=1,\dots,n,~~
 p_j(\tz)\equiv 1, \tz\in \mathcal{M},~j=1,2,3,
 \end{align}
 \begin{align}\label{gj}
 \E{\LT(W(\s)-W(\t)\RT)^2}=\sum_{i=1}^{n}|s_i-t_i|^{\alpha_i},~
 g_j(\bt_j)=\sum_{i\in \Lambda_j} b_it_i^{\beta_i},~j=1,2,3.
 \end{align}
 Let $\Gamma(x)=\int_{0}^{\infty} s^{x-1}e^{-s} d s$ for $x>0$ and for $\alpha\in (0,2]$ and $b>0$, define Pickands and Piterbarg constants respectively, for $\lambda>0$,
\begin{align}\label{H1}
\mathcal{H}_{B^{\alpha}}[0,\lambda] &=\E{\sup_{t\in [0,\lambda]}
	e^{ \sqrt{2}B^\alpha (t)-t^\alpha}},~~\mathcal{H}_{B^{\alpha}}=\lim_{\lambda\rw\IF}\frac{\mathcal{H}_{B^{\alpha}}[0,\lambda]}{\lambda},\nonumber\\
\mathcal{P}_{B^{\alpha}}^{b}[0,
\lambda]&=\E{\sup_{t\in [0,\lambda]}
	e^{ \sqrt{2}B^\alpha (t)-(1+b)t^\alpha}},~~
\mathcal{P}_{B^{\alpha}}^{b}=\lim_{\lambda\rw\IF}\mathcal{P}_{B^{\alpha}}^{b}[0,
\lambda],
\end{align}
where $B^{\alpha}$ represents a standard fractional Brownian motion with zero mean and covariance
$$Cov(B_{\alpha}(s),B_{\alpha}(t))=\frac{|t|^\alpha+|s|^{\alpha}-|t-s|^{\alpha}}{2},~s,t\geq 0.$$
We refer to \cite{Pit96} and the references therein for properties of Pickands and Piterbarg constants.

The following proposition partially generalizes Theorems 7.1 and 8.1 of \cite{Pit96}.

\begin{prop}\label{Prop1}
Under the assumption of Theorem \ref{MThm2}, if (\ref{ai})-(\ref{gj}) hold, then
\BQNY
\pk{\sup_{\t\in\mathcal{A}}X(\t)>u}\sim  C u^{\sum_{i\in \Lambda_0\cup\Lambda_1}\frac{2}{\alpha_i}-\sum_{i\in \Lambda_1}\frac{2}{\beta_i}}\Psi(u),
\EQNY
where
\BQNY
C=v_{k_0}(\mathcal{M})\left(\prod_{i\in \Lambda_0\cup\Lambda_1} a_i\mathcal{H}_{B^{\alpha_i}}\right) \left(\prod_{i\in\Lambda_1}b_i^{-1/\beta_i}\Gamma(1/\beta_i+1)\right)\prod_{i\in \Lambda_2} \mathcal{P}_{B^{\alpha_i}}^{a_i^{-\beta_i}b_i}.
\EQNY
\end{prop}
\section{Applications}
In this section, we illustrate our main results by application of Theorem \ref{MThm2}
to two classes of Gaussian fields
with nonstandard structures of the correlation function.

\subsection{The performance table, largest eigenvalue of GUE matrix and related problems}\label{s.GUE}
Let
\BQN\label{ZZ1}
Z^{\alpha}(\t):=\sum_{i=1}^{n+1}a_i\LT(B^{\alpha}_{i}(t_i)-B^{\alpha}_{i}(t_{i-1})\RT),\ \t=(t_1\ldot t_{n}),
\EQN
where
$t_0=0,t_{n+1}=1$
and $B^{\alpha}_{i},\ i=1\ldot n+1$ are mutually independent fractional Brownian motions with Hurst index $\alpha/2\in(0,1)$
and $a_i>0,\ i=1\ldot n+1$.
We are interested in the asymptotics of
\BQN\label{aim1}
\pk{D_n^\alpha>u}=\pk{\sup_{\t\in\mathcal{S}_{n}}Z^\alpha(\t)>u}
\EQN
for large $u$, where
$\mathcal{S}_{n}=\{\t\in \R^n: 0\leq t_1\leq\cdots\leq t_n\leq 1\}$.  Without loss of generality, we assume
that $\max_{i=1,\dots, n+1} a_i=1$.

Random variable $D_n^\alpha$ arises in many
problems that are important both in theoretical and applied probability. In particular it
is strongly related with the notion of {\it performance table}. More precisely,
following \cite{Bar01},
let
$\vk{w}=(w_{ij}), i,j\geq 1$
be a family of independent random values
indexed by the integer points of the first quarter of the plane.
A monotonous path $\pi$ from $(i,j)$ to $(i',j'), i\leq i'; j\leq j'; i,j,i',j'\in\N$ is a sequence $(i,j)=(i_0,j_0), (i_1,j_1)\ldot (i_l,j_l)=(i',j')$ of length $k=i'+j'-i-j+1$, such that all lattice steps $(i_k,j_k)\rw (i_{k+1},j_{k+1})$ are of size one and (consequently) go to the North or to the East. The weight $\vk{w}(\pi)$ of such a path is just the sum of all entries of the array $\vk{w}$ along the path.
We define performance table $l(i,j), i,j \in \N$ as the array of largest pathweights from $(1,1)$ to $(i,j)$, that is
\BQNY
l(i,j)=\max_{\pi \ \text{from}\ (1,1)\ \text{to}\ (i,j)} \vk{w}(\pi).
\EQNY
\kk{If $\Var(w_{ij})\equiv v>0$ and $\E{w_{ij}}\equiv e$ for all $i,j$,} then 
\BQNY
D_{n,k}:=\frac{l(n+1,k)-ke}{\sqrt{k v}}
\EQNY
converges in law as $k\rw\IF$ to 
$D_n^1$ with $a_i\equiv 1$; see \cite{Bar01}.
We refer to \cite{Bar01,GlW91,Sri93} and references therein for applications of
performance tables in queueing theory and in interacting particle systems.
Notably, as observed in \cite{Bar01}, if $a_i\equiv 1$
then $D^1_n$ has the same law as the largest eigenvalue of an
$n$-dimensional {\it Gaussian Unitary Ensamble} random matrix, see \cite{Oco02} for details
and further relations with non-colliding Brownian motions.

Denote
\BQN\label{mmd}
\mathcal{N}=\{i:a_i=1, i=1\ldot n+1\},\
\MNB=\{i: a_i<1, i=1\ldot n+1\},\
\mathfrak{m}=\sharp\mathcal{N}.
\EQN
For $k^*=\max\{i\in\mathcal{N}\}$ and  $\x=(x_1\ldot x_{k^*-1},x_{k^*+1}\ldot x_{n+1})$,
we define

\begin{align}\label{defineW}
W(\x)&=\frac{\sqrt{2}}{2}\sum_{i\in\mathcal{N}}
\LT(B_i(s_i(\x))-\widetilde{B}_i(s_{i-1}(\x))\RT)+\frac{\sqrt{2}}{2}\sum_{i\in\MNB}
a_i\LT(B_{i}(s_i(\x))-B_i(s_{i-1}(\x))\RT),
\end{align}
where $B_i, \widetilde{B}_i$ are independent standard Brownian motions and
\BQNY
s_i(\x)=\LT\{
\begin{array}{ll}
x_i,& \text{if}\ i\in \mathcal{N}\ \text{and}\ i<k^*,\\
\sum_{j=\max\{k\in\mathcal{N}:k<i\}}^{i}x_j,& \text{if}\ i\in\MNB\ \text{and}\ i<k^*,\\
\sum_{j=i+1}^{n+1}x_j,& \text{if}\ i\geq k^*,
\end{array}
\RT.
\EQNY
with the convention that $\max\emptyset=1$.

\cLa{For $\mathfrak{m}$ given in \eqref{mmd}} define
\BQN\label{PP}
\mathcal{H}_W=\lim_{\lambda\rw\IF}\frac{1}{\lambda^{\mathfrak{m}-1}}
\E{\sup_{\x\in[0,\lambda]^n}e^{ \sqrt{2}W(\x)-\LT(\underset{i\neq k^*}{\sum_{i=1}^{n+1}}x_i\RT)}}.
\EQN
It appears that for $\alpha=1$ and $\mathfrak{m}< n+1$
the field $Z^1$ satisfies {\bf A1} with
$W$ as given in (\ref{defineW}).
Notably, it has stationary increments with respect to coordinates  $\mathcal{N}$
while the increments of $W$ are not stationary with respect to coordinates  $\mathcal{N}^c$;
see (\ref{yrr}) in the proof of the following proposition.
Moreover, we have then
    $\Lambda_0=\mathcal{N}$, $\Lambda_1=\emptyset$, $\Lambda_2=\mathcal{N}^c$, $\Lambda_3=\emptyset$.

\begin{prop}\label{ThmM1} For $Z^{\alpha}$ defined in (\ref{ZZ1}),
we have, as $u\rw\IF$,
\BQNY
\pk{\sup_{\t\in\mathcal{S}_{n}}Z^\alpha(\t)>u}\sim
\LT\{
\begin{array}{ll}
	C u^{(\frac{2}{\alpha}-1)n}\Psi\LT(\frac{u}{\sigma_*}\RT), & \alpha\in(0,1),\\
	\cLa{\frac{1}{(\mathfrak{m}-1)!}}\mathcal{H}_Wu^{2(\mathfrak{m}-1)}\Psi(u),& \alpha=1,\\
	\mathfrak{m}\Psi(u), & \alpha\in(1,2),
\end{array}
\RT.
\EQNY
where $\sigma_*=\LT(\sum_{i=1}^{n+1}a_i^{\frac{2}{1-\alpha}}\RT)^{\frac{1-\alpha}{2}}
$ and
 $$
 C= \LT( \mathcal{H}_{B^{\alpha}}\RT)^n\LT(\prod_{i=1}^n\LT(a_i^2+a_{i+1}^2\RT)^{\frac{1}{\alpha}}\RT)2^{(1-\frac{1}{\alpha})n}\LT(\frac{\pi}{\alpha(1-\alpha)}\RT)^{\frac{n}{2}}\sigma_{*}^{\frac{-(\alpha-2)^2 n}{(1-\alpha)\alpha}}
\left(\sum_{j=1}^{n+1}\prod_{i\neq j}a_i^{\frac{2}{\alpha-1}}\right)^{-\frac{1}{2}}.
$$

\end{prop}

\begin{remark}\label{prop41}
i) \cLa{If $1\leq \mathfrak{m}\leq n$}, then
$1
 \le \mathcal{H}_W\leq
\COM{ \lim_{\lambda\rw\IF}\frac{1}{\lambda^{m-1}}
\E{\sup_{\x\in[0,\lambda]^n}e^{ \underset{i\neq k^*}{\sum_{i=1}^{n}}\sqrt{2n}B_i(x_i)
-\LT(n\underset{i\neq k^*}{\sum_{i=1}^{n}}x_i\RT)-\sum_{i\in\MNB}\frac{1-a_i^2}{2}x_i}}=}
n^{\mathfrak{m}-1}
\prod_{i\in\MNB}\LT(1+\frac{2n}{1-a_i^2}\RT).
$\\
\cLa{ii) If $\mathfrak{m}=n+1$, then $\mathcal{H}_W=1$.}
\end{remark}
The proof of Remark \ref{prop41} is postponed to Appendix.

\subsection{Chi processes}\label{s.chi}

Consider {\it chi} process $\chi$ generated by a  process $X$, that is let
 \begin{align}\label{chi}\chi(t):=\sqrt{\sum_{i=1}^{n}X_i^2(t)},~t\in [0,1],
 \end{align}
where $X_i$, $i=1,...,n$, are iid copies of $X$.
We suppose that $\{X(t), t\in[0,1]\}$  is a centered Gaussian process with a.s. continuous trajectories,
standard deviation  function
\begin{align}\label{X}
\sigma_X(t)=\frac{1}{1+bt^\alpha},~t\in [0,1],~\text{for}~b>0
\end{align}
and
correlation function
\begin{align}\label{Y}1-r(s,t)\sim aVar(Y(t)-Y(s)),~ s,t\to 0,~ \text{for}~ a>0,
\end{align}
where
$\{Y(t),~t\geq 0\}$ is a centered Gaussian process with a.s. continuous trajectories and satisfies:

{\bf B1}: $\{Y(t),~t\geq 0\}$ is  self-similar  with index $\alpha/2\in (0,1)$
(i.e. for all $r>0$,
$\{Y(rt),~t\geq 0\}\laweq\{r^{\alpha/2}Y(t),~t\geq 0\},$
where $\laweq$ means the equality of finite dimensional distributions) and $\sigma_Y(1)=1$;

{\bf B2}: there exist $c_Y>0$ and $\gamma \in [\alpha, 2]$ such that $$Var(Y(1)-Y(t))\sim c_Y|1-t|^\gamma,~t\uparrow 1.$$

Examples of Gaussian processes satisfying {\bf B1} and {\bf B2} cover such classes of Gaussian processes, as
fractional Brownian motions,
bi-fractional Brownian motions (see e.g. \cite{Houdre2003,Lei2009}), sub-fractional Brownian motions
(see e.g. \cite{Bojdecki2004,Dzhaparidze2004}), dual-fractional Brownian motions (see \cite{Li2004})
and time-average of fractional Brownian motions (see \cite{Li2004, Debic2020}).

For a Gaussian process $Y$ satisfying {\bf B1-B2}, a generalized Piterbarg constant is defined as, for $b>0$,
 \begin{align*}
 \mathcal{P}_{Y}^{b}=\lim_{S\to\infty}\E{\sup_{t\in [0,S]}e^{\sqrt{2}Y(t)-(1+b)t^{\alpha}}}\in (0,
 \infty).
 \end{align*}
We refer to \cite{Debic2020} for the finiteness and other properties of this constant.

The literature on the asymptotics of
\begin{eqnarray}
\pk{\sup_{t\in[0,1]}\chi(t)>u},\label{chi}
\end{eqnarray}
as $u\to\infty$
is focused on the case where $Y$ in (\ref{Y}) is a fractional Brownian motion, i.e.,
$1-r(s,t)\sim a|t-s|^\alpha$ as $s,t\to$ for some $\alpha \in (0,2]$,
which means that the correlation function of $X$ is locally homogeneous at $0$;
see \cite{Pitchi1994,EnkelejdJi2014Chi,LJ2017, Pit21}.
In the following proposition, $Y$ can be a general self-similar Gaussian  process satisfying {\bf B1-B2},
which allows for locally nonhomogeneous structures of the correlation function of $X$,
that were not investigated in the literature.

The idea of getting the asymptotics of (\ref{chi}) is based on
a transformation into supremum of Gaussian random field over a sphere, see \cite{Fatalov1993,Pitchi1994, Pit21}. That is,
we observe that
$$\sup_{t\in [0,1]}\chi(t)=\sup_{t\in [0,1], \sum_{i=1}^{n}v_i^2=1}X_i(t)v_i.$$
Next, we transform the Euclidean coordinates into spherical coordinates,
$$v_1(\vk{\theta})=\cos(\theta_1), v_2(\vk{\theta})=\sin(\theta_1)\cos(\theta_2), \dots, v_{n-1}(\vk{\theta})=\left(\prod_{i=1}^{n-2}\sin(\theta_i)\right)\cos(\theta_{n-1}), v_n(\vk{\theta})=\prod_{i=1}^{n-1}\sin(\theta_i),$$
where $\vk{\theta}=(\theta_1,\dots,\theta_{n-1})$ and $\vk{\theta}\in [0,\pi]^{n-2}\times [0,2\pi)$.
We denote
$$Z(\vk{\theta},t)=\sum_{i=1}^{n} X_i(t) v_i(\theta),~ \vk{\theta}\in [0,\pi]^{n-2}\times [0,2\pi), t\in [0,1]$$
and we have
$$\sup_{t\in [0,1]}\chi(t)=\sup_{(\vk{\theta}, t)\in E}Z(\vk{\theta},t)
~\text{with}~ E=[0,\pi]^{n-2}\times [0,2\pi)\times[0,1].$$
Consequently,
\begin{align}\label{chitrans}\mathbb{P}\left(\sup_{t\in [0,1]}\chi(t)>u\right)=\mathbb{P}\left(\sup_{(\vk{\theta}, t)\in E}Z(\vk{\theta},t)>u\right).
\end{align}
Then, it appears that the Gaussian field $Z$ satisfies Theorem \ref{MThm2} with
$W$ in (\ref{Rr}) and (\ref{W}) given by
\begin{align*}
W(\vk{\theta}, t)=\sum_{i=1}^{n-1}B_i^2(\theta_i)+\sqrt{a}Y(t),~\vk{\theta}\in \mathbb{R}^{n-1}\times \mathbb{R}^+,
\end{align*}
where $B_i^2$ are independent fractional Brownian motions with index $2$ and
$Y$ is the self-similar Gaussian process given in (\ref{Y})
that is independent of $B_i^2$.
Notably, one can check that, if $Y$ is not a fractional Brownian motion, the above defined
$W$ does not have stationary increments with respect to coordinate $t$.
Moreover, 
$\Lambda_0=\{1,\dots, n-1\}$, $\Lambda_1=\emptyset$, $\Lambda_2=\{n\}$, $\Lambda_3=\emptyset$.

 \begin{prop}\label{ex02}
 For $\chi$ defined in (\ref{chi}) with $X$ satisfying (\ref{X}) and (\ref{Y}), we have
\BQNY
\pk{\sup_{t\in[0,1]}\chi(t)>u}\sim \frac{2^{\frac{3-n}{2}}\sqrt{\pi}}{\Gamma(n/2)}\mathcal{P}_{Y}^{a^{-1}b}u^{n-1}\Psi(u),~u\to\infty.
\EQNY
\end{prop}

\def\MT{\mathcal{T}_{n}}
\def\MTK{\mathcal{T}_{n}^k}
\section{Proof of Theorem \ref{MThm2}}
We denote by $\mathbb{Q}, \mathbb{Q}_i, i=1,2,3,\dots$ positive constants that may differ from line to line.
\subsection{An  adapted  version of Theorem 2.1 in \cite{Uniform2016}}
In this subsection, we display a version of Theorem 2.1 in \cite{Uniform2016}, which plays an important role in the proof of Theorem \ref{MThm2}.
 Let $X_{u,\l}(\t), \t\in E\subset \mathbb{R}^n, \l\in K_u \subset \mathbb{R}^m,~ m\geq 1$   be a family of Gaussian random fields with variance $1$, where $E\subset \mathbb{R}^n$ is a compact set containing ${\bf 0}$ and $K_u\neq \emptyset$. Moreover, assume that $g_{u,\l}, \l\in K_u$ is a series of functions over $E$ and $u_{\l}, \l\in K_u$ are positive functions of $u$ satisfying $\lim_{u\rw\IF}\inf_{\l\in K_u}u_{\l}=\IF$.
In order to get the uniform asymptotics of
$$\pk{\sup_{\t\in E}\frac{X_{u,\l}(\t)}
 {1+g_{u,\l}(\t)}>u_{\l}}$$
 with respect to $\l\in K_u$,
 we shall impose the following assumptions:

{\bf C1}: There exists a function $g$ such that
\BQNY
\lim_{u\rw\IF}\sup_{\l\in K_u}\sup_{\t\in E}\left|u_{\l}^2g_{u,\l}(\t)-g(\t)\right|=0.
\EQNY
{\bf C2}: There exists a centered Gaussian random field $V(\t), \t\in E$ with $V({\bf 0})=0$ such that
\BQNY
\lim_{u\rw\IF}\sup_{\l\in K_u}\sup_{\s,\t\in E}\left|u_{\l}^2Var(X_{u,\l}(\t)-X_{u,\l}(\s))-2Var(V(\t)-V(\s))\right|=0.
\EQNY
{\bf C3}: There exist $\gamma\in (0,2]$  and $\mathcal{C}>0$ such that for $u$ sufficiently large
\BQNY
\sup_{\l\in K_u}\sup_{\s\neq \t, \s,\t\in E}u_{\l}^2\frac{Var(X_{u,\l}(\t)-X_{u,\l}(\s))}{\sum_{i=1}^n|s_i-t_i|^{\gamma}}\leq \mathcal{C}.
\EQNY

\BEL\label{uniform}
Let  $X_{u,\l}(\t), \t\in E\subset \mathbb{R}^n, \l\in K_u$  be a family of Gaussian random fields with variance $1$, $g_{u,\l}, \l\in K_u$ be functions defined on  $E$ and $u_{\l}, \l\in K_u$ be positive constants.  If {\bf C1-C3} are satisfied, then
\BQNY
\lim_{u\rw\IF}\sup_{\l\in K_u}\abs{\frac{\pk{\sup_{\t\in E}\frac{X_{u,\l}(\t)}
 {1+g_{u,l}(\t)}>u_{\l}}}{\Psi(u_{\l})}- \mathcal{P}^{g}_V\LT(E\RT)}=0,
\EQNY
where
\BQNY
\mathcal{P}^g_V\LT(E\RT)=
\E{\sup_{\t\in E}e^{ \sqrt{2}V(\t)-\sigma_V^2(\t)-g(\t)}}.
\EQNY
\EEL
\subsection{Proof of Theorem \ref{MThm2}}

 In order to simplify the proof, without loss of generality, we suppose that $\Lambda_0=\{1,\dots,k_0\}$ and  $\Lambda_i=\{k_{i-1}+1,\dots, k_i\}, ~i=1,2,3$, which in fact can be obtained by \cLa{change} of order of the coordinates.
\cLa{Thus} we have
 $\mathcal{M}^*=\{\t\in\mathcal{A}: t_{i}=0,i=k_0+1\ldot n\}$ and $\mathcal{M}=\{\tt:\t\in\mathcal{A}, t_{i}=0,i=k_0+1\ldot n\}$.
 The  proof is divided into three steps: In the first step, we show that the supremum  of $X(t)$ over $\mathcal{A}$ is predominately achieved on a subset; In the second step, we split this subset into small hyperrectangles and derive the asymptotics on each \cLa{hyperrectangle} resorting to the so-called double-sum method in \cite{Pit96}; Finally,  step 3 is devoted to adding up the asymptotics in step 2 to obtain the asymptotics over the whole set.

\subsubsection{\bf Step 1} In the first step of the proof, we divide $\mathcal{A}$ into two sets:
$$E_2(u)=\{\t\in\mathcal{A}:t_{i}\in[0,\delta_i(u)], k_0+1\leq i\leq n\}, \quad \delta_i(u)=\LT(\frac{\ln u}{u}\RT)^{2/\beta_i}, \quad k_0+1\leq i\leq n,$$
a neighborhood of $\mathcal{M}^*$, which \cLa{maximizes} the variance of $X(t)$ (with high probability the supremum is realized in $E_2(u)$) and the set $\mathcal{A}\setminus E_2(u)$, over which the probability associated with supremum is asymptotically negligible. For the lower bound, we only consider the process over
$$E_1(u)=\{\t\in\mathcal{A}:t_{i}\in[0,\delta_i(u)], k_0+1\leq i\leq k_1; t_{i}\in[0, u^{-2/\alpha_i}\lambda], k_1+1\leq i\leq k_2 ; t_{i}=0, k_2+1\leq i\leq k_3\}, ~  \lambda>0, $$  a neighborhood of $\mathcal{M}^*$.

To simplify  notation,  we denote for $\Delta_1, \Delta_2 \subseteq\R^{n}$
\BQNY
\pu{\Delta_1}:=\pk{\sup_{\t\in\Delta_1}X(\t)>u},\quad
\pu{\Delta_1,\Delta_2 }:=\pk{\sup_{\t\in\Delta_1}X(\t)>u,
\sup_{\t\in\Delta_2}X(\t)>u}.
\EQNY
Then  we have that  for any $u>0$
\BQN\label{boundlem1}
  \quad \pu{E_1(u)}\leq \pu{\mathcal{A}}\leq \pu{E_2(u)}+\pu{\mathcal{A}\setminus E_2(u)}.
\EQN

Note that in light of  \cite{Pit96} [Theorem 8.1], by (\ref{boundrr1}) and (\ref{G}),   for $u$ sufficiently large,
\BQN\label{err0}
\pu{\mathcal{A}\setminus E_2(u)}
\leq\mathbb{Q} v_n(\mathcal{A})u^{\sum_{i=1}^n\frac{2}{\alpha_i}}\Psi\LT(\frac{u}{1-\mathbb{Q}_1\LT(\frac{\ln u}{u}\RT)^2}\RT).
\EQN
\subsubsection{\bf Step 2}
In the second step, we divide $\mathcal{M}$ onto small hypercubes such that
$$\bigcup_{\vk{r}\in V^{-}}\mathcal{M}_{\vk{r}}\subset \mathcal{M}\subset \bigcup_{\vk{r}\in V^{+}}\mathcal{M}_{\vk{r}},$$
where
$$\mathcal{M}_{\vk{r}}=\prod_{i=1}^{k_0}[r_iv,(r_i+1)v], \quad \vk{r}=(r_1,\dots,r_{k_0}), r_i\in \mathbb{Z}, 1\leq i\leq k_0, v>0,$$
and
$$V^{+}:=\{\vk{r}: \mathcal{M}_{\vk{r}}\cap \mathcal{M}\neq \emptyset\}, \quad V^{-}:=\{\vk{r}: \mathcal{M}_{\vk{r}}\subset\mathcal{M}\}.$$
For fixed $\vk{r}$, we analyze the supremum of  $X$ over a set related to $\mathcal{M}_{\vk{r}}$. For this, let
$$E_{1,\vk{r}}(u)=\{\t: \tt\in \mathcal{M}_{\vk{r}}; t_{i}\in[0,\delta_i(u)], k_0+1\leq i\leq k_1; t_{i}\in[0, u^{-2/\alpha_i}\lambda], k_1+1\leq i\leq k_2 ; t_{i}=0, k_2+1\leq i\leq k_3\}, $$
$$E_{2,\vk{r}}(u)=\{\t:\tt\in \mathcal{M}_{\vk{r}}; t_{i}\in[0,\delta_i(u)], k_0+1\leq i\leq n\}.$$
Moreover, define an auxiliary set
$$
E_{3,\vk{r}}(u)=\{ (\tt, \bt_1,\bt_2):\tt\in\mathcal{M}_{\vk{r}}, t_{i}\in[0,\delta_i(u)], k_0+1\leq i\leq k_2 \}.$$
\cLa{We next focus on $\pu{E_{1,\vk{r}}(u)}$ and $\pu{E_{2,\vk{r}}(u)}$.
The idea of the proof of this step is first to split $E_{1,\vk{r}}(u)$ and $E_{2,\vk{r}}(u)$ onto tiny
hyperrectangles and uniformly derive the tail probability asymptotics on each hyperrectangle;
and then
to show that the asymptotics over $E_{i,\vk{r}}(u),i=1,2$ are the summation of the asymptotics over the corresponding
hyperrectangles, respectively.}

To this end, we   introduce the following notation. For some $\lambda>0$, let
\BQNY
I_{u,i}(l)=\LT[l\frac{\lambda}{u^{2/\alpha_i}},(l+1)\frac{\lambda}{u^{2/\alpha_i}}\RT],\ l\in\N, \quad \l=(l_1,\dots,l_n), \l_j=(l_{k_{j-1}+1},\dots, l_{k_j}),~ j=1,2,
\EQNY
\BQNY\Du(\l)=\left(\prod_{i=1}^{k_2}I_{u,i}(l_i)\right)\times \prod_{i=k_2+1}^n[0, \epsilon u^{-2/\alpha_i}],\quad
\mathcal{C}_u(\l)=\left(\prod_{i=1}^{k_1}I_{u,i}(l_i)\right)\times \prod_{i=k_1+1}^{k_2}[0, \lambda u^{-2/\alpha_i}]\times \overline{\vk{0}}_3,
\EQNY
with $\overline{\vk{0}}_3=(0,\dots,0)\in\mathbb{R}^{n-k_2}$, and
\BQNY
M_i(u)=\LT\lfloor\frac{vu^{2/\alpha_i}}{\lambda}\RT\rfloor, 1\leq i\leq k_0, \quad M_i(u)=\LT\lfloor\frac{\delta_i(u)u^{2/\alpha_i}}{\lambda}\RT\rfloor, k_0+1\leq i\leq k_2.
\EQNY

\COM{In view of  \eqref{lem2rr1} and {\bf B1}, we have
\BQN\label{rr111}
\lim_{\delta\rw 0}\underset{\abs{\s-\t}\leq \delta}{\sup_{\s,\t\in \mathcal{A}}}\abs{\frac{1-r(\s,\t)}{\sigma^2_W(\abs{t_1-s_1}\ldot \abs{t_n-s_n})}-1}=0.
\EQN}
In order to derive an upper bound for $\pu{E_{2,\vk{r}}(u)}$ and a lower bound for $\pu{E_{1,\vk{r}}(u)}$, we introduce the following notation
 for some $\epsilon\in(0,1)$,
\BQNY
\mathcal{L}_{1}(u)&=&\LT\{\l:\prod_{i=1}^{k_2}I_{u,i}(l_i)\subset E_{3,\vk{r}}(u), l_i=0, k_1+1\leq i\leq n\RT\},\\
\mathcal{L}_{2}(u)&=&\LT\{\l:\left(\prod_{i=1}^{k_2}I_{u,i}(l_i)\right)\cap E_{3,\vk{r}}(u)\neq\emptyset, l_i=0, k_1+1\leq i\leq n\RT\},\\
\mathcal{L}_{3}(u)&=&\LT\{\l:\left(\prod_{i=1}^{k_2}I_{u,i}(l_i)\right)\cap E_{3,\vk{r}}(u)\neq\emptyset, \sum_{i=k_1+1}^{k_2}l_i^2>0, l_i=0, k_2+1\leq i\leq n\RT\},\\
\mathcal{K}_1(u)&=&\LT\{(\l,\j):\l,\j\in\mathcal{L}_1(u),\mathcal{C}_u(\l)\cap \mathcal{C}_u(\j)\neq\emptyset \RT\},\\
\mathcal{K}_2(u)&=&\LT\{(\l,\j):\l,\j\in\mathcal{L}_1(u),\mathcal{C}_u(\l)\cap \mathcal{C}_u(\j)=\emptyset \RT\},\\
u^{-\epsilon}_{\l_1}&=&u\LT(1+(1-\epsilon)\inf_{\bt_1\in [\l_1,\l_1+1]}p_{1,\vk{r}}^-g_1(u^{-2/\aalpha_1}\lambda\bt_1)\RT),\\
u^{+\epsilon}_{\l_1}&=&u\LT(1+(1+\epsilon)\sup_{\bt_1\in [\l_1,\l_1+1]}p_{1,\vk{r}}^+g_1(u^{-2/\aalpha_1}\lambda\bt_1)\RT),\\
p_{j,\vk{r}}^+&=&\sup_{\tz\in \mathcal{M}_{\vk{r}}}p_j(\tz), \quad p_{j,\vk{r}}^-=\inf_{\tz\in \mathcal{M}_{\vk{r}}}p_j(\tz), \quad j=1,2,3.
\EQNY
Bonferroni inequality gives that  for $u$ sufficiently large
\BQN\label{boundupperlower112}
\pu{E_{1,\vk{r}}(u)}&\geq &\sum_{\l\in\mathcal{L}_{1}(u)}\pu{\mathcal{C}_u(\l)}-\sum_{i=1}^{2}\Gamma_i(u),
\EQN
\BQN\label{boundupperlower111}
\pu{E_{2,\vk{r}}(u)}&\leq& \sum_{\l\in\mathcal{L}_{2}(u)}\pu{\Du(\l)}+\sum_{\l\in\mathcal{L}_{3}(u)}\pu{\Du(\l)},
\EQN
where
\BQNY
\Gamma_i(u)=\sum_{(\l,\j)\in\mathcal{K}_i(u)}\pu{\mathcal{C}_u(\l),\mathcal{C}_u(\j)},\ i=1,2.
\EQNY
We first derive the upper bound of $\pu{E_{2,\vk{r}}(u)}$ as $u\to\infty$. To this end, we need to find the upper bounds of $\sum_{\l\in\mathcal{L}_{j}(u)}\pu{\Du(\l)}, j=2,3$, separately.

{\it \underline{Upper bound for $\sum_{\l\in\mathcal{L}_{2}(u)}\pu{\Du(\l)}$}}.
By \eqref{Var}, we have that for $u$ sufficiently large
\BQNY
\sum_{\l\in\mathcal{L}_{2}(u)}\pu{\Du(\l)}
& \leq&
 \sum_{\l\in\mathcal{L}_{2}(u)}\pk{\sup_{\t\in\Du(\l)}\frac{\overline{X}(\t)}{1+(1-\epsilon)p_{2,\vk{r}}^-g_2(\bt_2)}
 >u_{\l_1}^{-\epsilon}}\\
 & =&\sum_{\l\in\mathcal{L}_{2}(u)}\pk{\sup_{\t\in E(\l,u)}\frac{X_{u,\l}(\t)}{1+(1-\epsilon)p_{2,\vk{r}}^-g_2(u^{-2/\aalpha_2}(\ta_2(\tz(\l,u)))^{-1}\bt_2)}
 >u_{\l_1}^{-\epsilon}},
\EQNY
where
\BQN\label{Xul}
X_{u,\l}(\t)=\overline{X}\left(u^{-2/\alpha_1}(l_1\lambda+(a_{1}(\tz(\l,u)))^{-1}t_1),\dots,  u^{-2/\alpha_n}(l_n\lambda+(a_{n}(\tz(\l,u))^{-1}t_n)\right),
\EQN
with $\tz(\l,u)=(u^{-2/\alpha_1}l_1, \dots, u^{-2/\alpha_{k}}l_k)$
and $E(\l,u)=\left(\prod_{i=1}^{k_2}[0, a_i(\tz(\l,u))\lambda]\right)\times \prod_{i=k_2+1}^{n}[0, a_i(\tz(\l,u))\epsilon].$

Note that by (\ref{G}),
$$u^{-2}g_{2,\vk{r}}^-(\bt_2)\leq g_2(u^{-2/\aalpha_2}(\ta_2(\tz(\l,u)))^{-1}\bt_2)=u^{-2}g_2((\ta_2(\tz(\l,u)))^{-1}\bt_2)\leq u^{-2}g_{2,\vk{r}}^+(\bt_2),$$
where
$$g_{2,\vk{r}}^-(\bt_2)=\inf_{\tz\in \mathcal{M}_{\vk{r}}}g_2((\ta_2(\tz)^{-1}\bt_2), \quad g_{2,\vk{r}}^+(\bt_2)=\sup_{\tz\in \mathcal{M}_{\vk{r}}}g_2((\ta_2(\tz)^{-1}\bt_2).$$
Moreover,
$$  E_{\vk{r}}^-\subset E(\l,u)\subset E_{\vk{r}}^+,$$
where
$$E_{\vk{r}}^+:=\left(\prod_{i=1}^{k_2}[0, a_{i,\vk{r}}^+\lambda]\right)\times \prod_{i=k_2+1}^{n}[0,a_{i,\vk{r}}^+\epsilon], \quad E_{\vk{r}}^-:=\left(\prod_{i=1}^{k_2}[0, a_{i,\vk{r}}^-\lambda]\right)\times \prod_{i=k_2+1}^{n}[0,a_{i,\vk{r}}^-\epsilon]$$
with
$$a_{i,\vk{r}}^+=\sup_{\tz\in \mathcal{M}_{\vk{r}}}a_i(\tz), \quad a_{i,\vk{r}}^-=\inf_{\tz\in \mathcal{M}_{\vk{r}}}a_i(\tz).$$
Hence
\BQN\label{sum1}
\sum_{\l\in\mathcal{L}_{2}(u)}\pu{\Du(\l)}\leq \sum_{\l\in\mathcal{L}_{2}(u)}\pk{\sup_{\t\in E_{\vk{r}}^+}\frac{X_{u,\l}(\t)}{1+(1-\epsilon)u^{-2}p_{2,\vk{r}}^-g_{2,\vk{r}}^-(\bt_2)}
 >u_{\l_1}^{-\epsilon}}.
\EQN

{\it \underline{Uniform asymptotics for the summands in (\ref{sum1})}}.  We need to specify the notation in Lemma \ref{uniform} for the current case.  \cLa{Let $X_{u,\l}$ be as was defined in (\ref{Xul}) and let}
$$u_{\l}=u_{\l_1}^{-\epsilon},\quad g_{u,\l}(\t)=(1-\epsilon)u^{-2}p_{2,\vk{r}}^-g_{2,\vk{r}}^-(\bt_2),\quad K_{u}=\mathcal{L}_{2}(u).$$
We first note that
$\lim_{u\rw\IF}\inf_{\l\in \mathcal{L}_2(u)}u_{\l_1}^{-\epsilon}=\IF$,
which  combined with continuity of $g_2$ \cLa{implies}
$$\lim_{u\rw\IF}\sup_{\l\in K_u}\sup_{\t\in E_{\vk{r}}^+}\left|u_{\l}^2g_{u,\l}(\t)-(1-\epsilon)p_{2,\vk{r}}^-g_{2,\vk{r}}^-(\bt_2)\right|=0.$$
Hence {\bf C1} holds with $g(\bt)=(1-\epsilon)p_{2,\vk{r}}^-g_{2,\vk{r}}^-(\bt_2)$.
To check {\bf C2}, by (\ref{Rr}) and (\ref{W}) and using the homogeneity of the increments of $W$ for fixed $\bt_2$ and $ \bt_3$, we have
\BQNY
\lim_{u\rw\IF}\sup_{\l\in K_u}\sup_{\s,\t\in E_{\vk{r}}^+}\left|u_{\l}^2Var(X_{u,\l}(\t)-X_{u,\l}(\s))-2Var(W(\t)-W(\s))\right|=0.
\EQNY
This implies that {\bf C2} is satisfied with the limiting stochastic process $W$ defined in {\bf A1}. {\bf C3} follows directly from (\ref{boundrr1}). Therefore, we conclude that
\BQN\label{uniform1}
\lim_{u\rw\IF}\sup_{\l\in K_u}\left|\frac{\pk{\sup_{\t\in E_{\vk{r}}^+}\frac{X_{u,\l}(\t)}{1+(1-\epsilon)u^{-2}p_{2,\vk{r}}^-g_{2,\vk{r}}^-(\bt_2)}
 >u_{\l_1}^{-\epsilon}}}{\Psi(u_{\l}^{-\epsilon})}-\mathcal{H}_W^{(1-\epsilon)p_{2,\vk{r}}^-g_{2,\vk{r}}^-(\bt_2)}\LT(E_{\vk{r}}^+\RT)\right|=0,
\EQN
where
\BQN\label{HW}
\mathcal{H}_W^{(1-\epsilon)p_{2,\vk{r}}^-g_{2,\vk{r}}^-(\bt_2)}(E_{\vk{r}}^+)=\E{\sup_{\t\in E_{\vk{r}}^+}e^{ \sqrt{2}W(\t)-\sigma^2_W(\t)-(1-\epsilon)p_{2,\vk{r}}^-g_{2,\vk{r}}^-(\bt_2)}}.
\EQN
Hence we have
\BQN\label{lem1bound111}
&&\sum_{\l\in\mathcal{L}_{2}(u)}\pk{\sup_{\t\in E}X_{u,\l}(\t)
 >u_{\l}^{-\epsilon}}\nonumber\\
&&\quad\leq
 \sum_{\l\in\mathcal{L}_{2}(u)} \mathcal{H}_W^{(1-\epsilon)p_{2,\vk{r}}^-g_{2,\vk{r}}^-(\bt_2)}(E_{\vk{r}}^+)  \Psi(u_{\l}^{-\epsilon})\nonumber\\
 &&\quad \leq \mathcal{H}_W^{(1-\epsilon)p_{2,\vk{r}}^-g_{2,\vk{r}}^-(\bt_2)}(E_{\vk{r}}^+) \Psi(u)\LT(\prod_{i=1}^{k_0}\frac{vu^{2/\alpha_i}}{\lambda}\RT)\sum_{i=k_0+1}^{k_1}\sum_{l_i=0}^{M_i(u)} e^{-(1-\epsilon)\inf_{\bt_1\in [\l_1,\l_1+1]}p_{1,\vk{r}}^-g_1(u^{2/\bbeta_1-2/\aalpha_1}\lambda\bt_1)}\nonumber\\
&&\quad \sim\frac{\mathcal{H}_W^{(1-\epsilon)p_{2,\vk{r}}^-g_{2,\vk{r}}^-(\bt_2)}(E_{\vk{r}}^+)}{\lambda^{k_1}} v^{k_0}\Psi(u)u^{\sum_{i=1}^{k_1}\frac{2}{\alpha_i}-\sum_{i=k_0+1}^{k_1}\frac{2}
{\beta_i}}\int_{\bt_1\in[0,\IF)^{k_1-k_0}} e^{-(1-\epsilon)p_{1,\vk{r}}^-g_1(\bt_1)}d\bt_1, \quad u\rw\IF.
\EQN

Note that
\BQNY\lim_{\epsilon\rw 0}\mathcal{H}_W^{(1-\epsilon)p_{2,\vk{r}}^-g_{2,\vk{r}}^-(\bt_2)}(E_{\vk{r}}^+)&=&\E{\sup_{(\tilde{\t}, \bt_1,\bt_2)\in\prod_{i=1}^{k_2}[0, a_{i,\vk{r}}^+\lambda]}e^{ \sqrt{2}W(\tilde{\t}, \bt_1,\bt_2,\bar{\vk{0}}_3)-\sigma^2_W(\tilde{\t}, \bt_1,\bt_2,\bar{\vk{0}}_3)-p_{2,\vk{r}}^-g_{2,\vk{r}}^-(\bt_2)}}\\
&:=&\mathcal{H}_W^{p_{2,\vk{r}}^-g_{2,\vk{r}}^-(\bt_2)}\left(\prod_{i=1}^{k_2}[0, a_{i,\vk{r}}^+\lambda]\right)
\EQNY
and by dominated convergence theorem, it follows that
$$\lim_{\epsilon\rw 0}\int_{\bt_1\in[0,\IF)^{k_1-k_0}} e^{-(1-\epsilon)p_{1,\vk{r}}^-g_1(\bt)}d\bt_1=\int_{\bt_1\in[0,\IF)^{k_1-k_0}} e^{-p_{1,\vk{r}}^-g_1(\bt)}d\bt_1.$$
\cLa{Hence}, letting $\epsilon\rw 0$ in (\ref{lem1bound111}), we have
\BQN\label{upper1}
\sum_{\l\in\mathcal{L}_{2}(u)}\pu{\Du(\l)}
 \leq \frac{\mathcal{H}_W^{p_{2,\vk{r}}^-g_{2,\vk{r}}^-(\bt_2)}\left(\prod_{i=1}^{k_2}[0, a_{i,\vk{r}}^+\lambda]\right)}{\lambda^{k_1}}v^{k_0}\Theta^-(u),\quad u\rw\IF,
\EQN
where
$$\Theta^\pm(u):=\Psi(u)u^{\sum_{i=1}^{k_1}\frac{2}{\alpha_i}-\sum_{i=k+1}^{k_1}\frac{2}{\beta_i}}\int_{\bt_1\in[0,\IF)^{k_1-k}} e^{-p_{1,\vk{r}}^\pm g_1(\bt)}d\bt_1.$$
{\it \underline{Upper bound for $\sum_{\l\in\mathcal{L}_{3}(u)}\pu{\Du(\l)}$}}.
\cLa{Next we find a tight asymptotic upper bound for the second term displayed in the right hand side of (\ref{boundupperlower111})}. For $u$ sufficiently large
\BQNY
\sum_{\l\in\mathcal{L}_{3}(u)}\pu{\Du(\l)}
&\leq&
 \sum_{\l\in\mathcal{L}_{3}(u)}\pk{\sup_{\t\in\Du(\l)}\overline{X}(\t)
 >u_{\l_1,\l_2}^{-\epsilon}}
 =\sum_{\l\in\mathcal{L}_{3}(u)}\pk{\sup_{\t\in E}\widetilde{X}_{u,\l}(\t)
 >u_{\l_1,\l_2}^{-\epsilon}},
\EQNY
where
\BQN\label{Xull}
\widetilde{X}_{u,\l}(\t)=\overline{X}\left(u^{-2/\alpha_1}(l_1\lambda+t_1),\dots,  u^{-2/\alpha_n}(l_n\lambda+t_n)\right),\quad E=[0,\lambda]^{k_2}\times [0,\epsilon]^{n-k_2},
\EQN
$$u_{\l_1,\l_2}^{-\epsilon}=u\left(1+(1-\epsilon)\inf_{\bt_1\in [\l_1, \l_1+1]}g_1(u^{-2/\aalpha_1}\lambda\bt_1)+(1-\epsilon)\inf_{\bt_2\in [\l_2, \l_2+1]}g_2(u^{-2/\aalpha_2}\lambda\bt_2)\right).$$
Let $Z_u(\t)$ be a homogeneous Gaussian random field with variance $1$ and correlation function satisfying
\BQN\label{Zu}
r_u(\s,\t)=e^{-u^{-2}2\mathcal{Q}_2\sum_{i=1}^n \abs{s_i-t_i}^{\alpha_i}}.
\EQN
By (\ref{boundrr1}) and Slepian's inequality (see, e.g., Theorem 2.2.1 in \cite{AdlerTaylor}) we have that for $u$ sufficiently large
$$\pk{\sup_{\t\in E}\widetilde{X}_{u,\l}(\t)
 >u_{\l_1,\l_2}^{-\epsilon}}\leq \pk{\sup_{\t\in E}Z_u(\t)
 >u_{\l_1,\l_2}^{-\epsilon}}, \quad \l\in \mathcal{L}_3(u).$$
Similarly as in the proof of (\ref{uniform1}), we have
\BQN\label{uniform2}
\lim_{u\rw\IF}\sup_{\l\in \mathcal{L}_3(u)}\left|\frac{\pk{\sup_{\t\in E}Z_u(\t)>u_{\l_1,\l_2}^{-\epsilon}}}{\Psi(u_{\l_1,\l_2}^{-\epsilon})}-\mathcal{J}(E)
\right|=0,
\EQN
where
\BQNY
\mathcal{J}(E)=\left(\prod_{i=1}^{k_2}\mathcal{H}_{B^{\alpha_i}}[0,(2\mathcal{Q}_2)^{1/\alpha_i}\lambda]\right)\left(\prod_{i=k_2+1}^n
\mathcal{H}_{B^{\alpha_i}}[0,\epsilon(2\mathcal{Q}_2)^{1/\alpha_i}\lambda]\right).
\EQNY
 Hence using the above asymptotics and assumption (\ref{G})
\BQNY
\sum_{\l\in\mathcal{L}_{3}(u)}\pu{\Du(\l)}&\leq& \sum_{\l\in\mathcal{L}_{3}(u)}\mathcal{J}(E)\Psi(u_{\l_1.\l_2}^{-\epsilon})\\
&\leq& \mathcal{J}(E)\Psi(u)\sum_{\l\in\mathcal{L}_{3}(u)}e^{-(1-\epsilon)\inf_{\bt_1\in [\l_1, \l_1+1]}u^2g_1(u^{-2/\aalpha_1}\lambda\bt_1)-(1-\epsilon)\inf_{\bt_2\in [\l_2, \l_2+1]}u^2g_2(u^{-2/\aalpha_2}\lambda\bt_2)}\\
&\leq&\mathcal{J}(E)\Psi(u)\LT(\prod_{i=1}^{k_0}\frac{vu^{2/\alpha_i}}
{\lambda}\RT)\sum_{i=k_0+1}^{k_1}\sum_{l_i=0}^{M_i(u)} e^{-(1-\epsilon)\inf_{\bt_1\in [\l_1, \l_1+1]}g_1(u^{2/\bbeta_1-2/\aalpha_1}\lambda\bt_1)}\\
&&\quad \times \sum_{l_{k_1+1}^2+\dots+l_{k_2}^2\geq 1, l_i\geq 0, k_1+1\leq i\leq k_2}e^{-(1-\epsilon)\inf_{\bt_2\in [\l_2, \l_2+1]}g_2(u^{2/\bbeta_2-2/\aalpha_2}\lambda\bt_2)}.
\EQNY
Moreover,
direct calculation shows that
\BQNY
\sum_{i=k_0+1}^{k_1}\sum_{l_i=0}^{M_i(u)} e^{-(1-\epsilon)\inf_{\bt_1\in [\l_1, \l_1+1]}g_1(u^{2/\bbeta_1-2/\aalpha_1}\lambda\bt_1)}\sim u^{\sum_{i=k_0+1}^{k_1}(\frac{2}{\alpha_i}-\frac{2}{\beta_i})}\lambda^{k_0-k_1}\int_{\bt_1\in[0,\IF)^{k_1-k_0}} e^{-(1-\epsilon)g_1(\bt)}d\bt_1, \quad u\rw\IF.
\EQNY
By assumption (\ref{G}) and the fact that $\aalpha_2=\bbeta_2$, we have for $\lambda>1$
\BQNY
&&\sum_{l_{k_1+1}^2+\dots+l_{k_2}^2\geq 1, l_i\geq 0, k_1+1\leq i\leq k_2}e^{-(1-\epsilon)\inf_{\bt_2\in [\l_2, \l_2+1]}g_2(u^{2/\bbeta_2-2/\aalpha_2}\lambda\bt_2)}\\
&&\leq \sum_{l_{k_1+1}^2+\dots+l_{k_2}^2\geq 1, l_i\geq 0, k_1+1\leq i\leq k_2}e^{-(1-\epsilon)c_{2,1}\sum_{i=k_1+1}^{k_2}(l_i\lambda)^{\beta_i}}\leq \mathbb{Q}_3e^{-\mathbb{Q}_2\lambda^{\beta^*}},
\EQNY
where $\beta^*=\min_{i=k_1+1}^{k_2}(\beta_i).$ Additionally,
\BQNY
\lim_{\epsilon\rw 0}\mathcal{J}(E)=\prod_{i=1}^{k_2}\mathcal{H}_{B^{\alpha_i}}[0,(2\mathcal{Q}_2)^{1/\alpha_i}\lambda],
\EQNY
and for $\lambda>1$
\BQNY
\prod_{i=1}^{k_2}\mathcal{H}_{B^{\alpha_i}}[0,(2\mathcal{Q}_2)^{1/\alpha_i}\lambda]\leq \mathbb{Q}_3\lambda^{k_2}.
\EQNY
Thus we have that for $\lambda>1$
\BQN\label{upper2}
\sum_{\l\in\mathcal{L}_{3}(u)}\pu{\Du(\l)}
 \leq \mathbb{Q}_3\lambda^{k_2-k_1}e^{-\mathbb{Q}_2\lambda^{\beta^*}}v^{k_0}\Theta^-(u),\quad u\rw\IF.
\EQN
{\it \underline{Upper bound for $\pu{E_{2,\vk{r}}(u)}$}}.
Combination of (\ref{upper1}) and (\ref{upper2}) yields that for $\lambda>1$
\BQN\label{upper}
\pu{E_{2,\vk{r}}(u)}\leq \left(\frac{\mathcal{H}_W^{p_{2,\vk{r}}^-g_{2,\vk{r}}^-(\bt_2)}\left(\prod_{i=1}^{k_2}[0, a_{i,\vk{r}}^+\lambda]\right)}{\lambda^{k_1}}
+\mathbb{Q}_3\lambda^{k_2-k_1}e^{-\mathbb{Q}_2\lambda^{\beta^*}}\right)v^{k_0}\Theta^-(u), \quad u\rw\IF.
\EQN
We next find a lower bound of $\pu{E_{1,\vk{r}}(u)}$ as $u\to\infty$ for which we need to derive the lower bound of $\sum_{\l\in\mathcal{L}_{1}(u)}\pu{\mathcal{C}_u(\l)}$  and upper bounds of $\Gamma_i(u), i=1,2$, respectively.

{\it \underline{Lower bound for $\sum_{\l\in\mathcal{L}_{1}(u)}\pu{\mathcal{C}_u(\l)}$}}.
\COM{In the derivation of the lower bound,  for simplicity, we write $X(\t)=X(\tt, \bt_1,\bt_2, \overline{\vk{0}}_3)$, \cLa{where $\overline{\vk{0}}_3=(0,\dots, 0)\in \mathbb{R}^{n-k_2}$}.}
Analogously to (\ref{upper1}), we derive
\BQN\label{lem1bound2111}
\sum_{\l\in\mathcal{L}_{1}(u)}\pu{\mathcal{C}_u(\l)}&\geq&
\frac{\mathcal{H}_W^{p_{2,\vk{r}}^+g_{2,\vk{r}}^+(\bt_2)}\left(\prod_{i=1}^{k_2}[0, a_{i,\vk{r}}^-\lambda]\right)}{\lambda^{k_1}}v^{k_0}\Theta^+(u),\quad u\rw\IF, \epsilon\rw 0.
\EQN

{\underline{Upper bound for $\Gamma_i(u), i=1,2$}}.
Applying \kk{an approach analogous} to that of the proof of Theorem 8.2 in \cite{Pit96}, we have, for $\lambda>1$ and as $u\to\infty$,
\BQN
\Gamma_1(u)
&\leq& \mathbb{Q}_4\lambda^{-1/2}
\lambda^{2k_2-k_1}
 v^{k_0}\Theta^-(u),\label{lem1bound4111}\\
\Gamma_{2}(u)
&\leq&\mathbb{Q}_{5}\lambda^{2k_2-k_1}e^{-\mathbb{Q}_{6}\lambda^{\alpha^*}}v^{k_0}\Theta^-(u),\label{lem1bound5111}
\EQN
where
$\alpha^*=\max(\alpha_1,\dots, \alpha_{k_1})$ and $\mathbb{Q}_i, i=4,5,6$ are some positive constants.

{\it \underline{Lower bound for $\pu{E_{1,\vk{r}}(u)}$}}.
Inserting  \eqref{lem1bound2111}, \eqref{lem1bound4111} and \eqref{lem1bound5111} into \eqref{boundupperlower112}, we obtain for $\lambda>1$, as $u\rw\IF$,
\BQN\label{lower}
\pu{E_{1,\vk{r}}(u)}\geq \left(\frac{\mathcal{H}_W^{p_{2,\vk{r}}^+g_{2,\vk{r}}^+(\bt_2)}\left(\prod_{i=1}^{k_2}[0, a_{i,\vk{r}}^-\lambda]\right)}{\lambda^{k_1}}-\mathbb{Q}_4\lambda^{-1/2}-\mathbb{Q}_{5}\lambda^{2k_2-k_1}
e^{-\mathbb{Q}_{6}\lambda^{\alpha^*}}\right) v^{k_0}\Theta^+(u).
\EQN
{\it\underline{Existence of $\mathcal{H}_W^{g_{2}(\bt_2)}$}}.
The idea used here is similar to that of Lemmas 7.1 and 8.3 in \cite{Pit96}.
\kk{Thus we present only main steps of argumentation}.  We
assume $$\vk{a}(\tz)=1,  p_j(\tz)=1, j=1,2,3, \tz\in \mathcal{M}_{\vk{r}}.$$
Dividing (\ref{upper}) and (\ref{lower}) by $v^{k_0}\Theta^-(u)$ and letting $u\rw\IF$, we derive that, for some $\lambda, \lambda_1>1$,
$$\limsup_{\lambda\rw\IF}\frac{\mathcal{H}_W^{g_{2}(\bt_2)}([0,\lambda]^{k_2})}{\lambda^{k_1}}
\leq\liminf_{\lambda\rw\IF}\frac{\mathcal{H}_W^{g_{2}(\bt_2)}([0,\lambda]^{k_2})}{\lambda^{k_1}}<\IF.$$
The positivity of the above limit follows from the same arguments as in \cite{Pit96}.
Therefore,
\BQN\label{Pickands}
\mathcal{H}_W^{g_{2}(\bt_2)}:=\lim_{\lambda\rw\IF} \frac{\mathcal{H}_W^{g_{2}(\bt_2)}([0,\lambda]^{k_2})}{\lambda^{k_1}}\in (0,
\IF).
\EQN
Moreover, using (\ref{upper}) and (\ref{lower}), we have, for $\lambda>1$,
\BQN\label{uniform bound}
\left|\frac{\mathcal{H}_W^{g_{2}(\bt_2)}([0,\lambda]^{k_2})}{\lambda^{k_1}}-\mathcal{H}_W^{g_{2}(\bt_2)}\right|
 \leq \mathbb{Q}_7\left(\lambda^{-1/2}+
\lambda^{2k_2-k_1}e^{-\mathbb{Q}_{6}\lambda^{\alpha^*}}
+\lambda^{k_2-k_1}e^{-\mathbb{Q}_2\lambda^{\beta^*}}\right).
\EQN

Let
$\mathcal{G}:=\{g_2: \text{$g_2$ is continuous}, ~ u g_2(\bt_2)=g_2(u^{1/{ \bbeta_{2}}}\bt_{2}), u>0,~ \inf_{\sum_{i=k_{1}+1}^{k_2}|t_i|^{\beta_i}=1}g_2(\bt_2)>c>0\},$
where $c$ and $\bbeta_{2}$ are fixed. For any $g_2\in\mathcal{G}$, (\ref{upper1}) and (\ref{upper2})-(\ref{lower}) are still valid. Hence, (\ref{uniform bound}) also holds. This implies that for any $\lambda>1$,
\BQN\label{uniform bound1}
\sup_{g_2\in \mathcal{G}}\left|\frac{\mathcal{H}_W^{g_{2}(\bt_2)}([0,\lambda]^{k_2})}{\lambda^{k_1}}-\mathcal{H}_W^{g_{2}(\bt_2)}\right|
 \leq \mathbb{Q}_7\left(\lambda^{-1/2}+
\lambda^{2k_2-k_1}e^{-\mathbb{Q}_{6}\lambda^{\alpha^*}}
+\lambda^{k_2-k_1}e^{-\mathbb{Q}_2\lambda^{\beta^*}}\right).
\EQN
\subsubsection{\bf Step 3} In this step of the proof, we sum up the asymptotics derived in step 2.
 Set
$$\Theta_1(u)=u^{\sum_{i=1}^{k_1}\frac{2}{\alpha_i}-\sum_{i=k+1}^{k_1}\frac{2}{\beta_i}}\Psi(u).$$
Letting $\lambda\rw \IF$ in (\ref{upper}) and (\ref{lower}), it follows that
\BQNY
\pu{E_{1,\vk{r}}(u)}
 &\geq& \mathcal{H}_W^{p_{2,\vk{r}}^+g_{2,\vk{r}}^+(\bt_2)}\prod_{i=1}^{k_1} a_{i,\vk{r}}^-\int_{\bt_1\in[0,\IF)^{k_1-k}} e^{-p_{1,\vk{r}}^+ g_1(\bt)}d\bt_1v^{k_0}\Theta_1(u),\\ \pu{E_{2,\vk{r}}(u)}
 &\leq& \mathcal{H}_W^{p_{2,\vk{r}}^-g_{2,\vk{r}}^-(\bt_2)}\prod_{i=1}^{k_1} a_{i,\vk{r}}^+\int_{\bt_1\in[0,\IF)^{k_1-k}} e^{-p_{1,\vk{r}}^- g_1(\bt)}d\bt_1v^{k_0}\Theta_1(u).
\EQNY

 We add up $\pu{E_{1,\vk{r}}(u)}$ and $\pu{E_{2,\vk{r}}(u)}$ with respect to $\vk{r}$ respectively to get a lower bound of $\pu{E_1(u)}$ and an upper bound of $\pu{E_2(u)}$. Observe that
\begin{align}
\pu{E_1(u)}&\geq \sum_{\vk{r}\in V^-}\pu{E_{1,\vk{r}}(u)}-\sum_{\vk{r}, \vk{r}'\in V^-, \vk{r}\neq \vk{r}'}\pu{E_{1,\vk{r}}(u), E_{1,\vk{r}'}(u)},\label{double}\\
\pu{E_2(u)}&\leq \sum_{\vk{r}\in V^+}\pu{E_{2,\vk{r}}(u)}.\nonumber
\end{align}

Note that $g_{2,\vk{r}}^+(\bt_2)\in \mathcal{G}, \vk{r}\in V^+$ and $p_2(\tz)g_2(\vk{a}_2^{-1}(\tz)\bt_2)\in \mathcal{G},\tz\in \mathcal{M}$ with  fixed $c$ and $\bbeta_{2}$.  Thus
(\ref{uniform bound1}) implies that for any $\epsilon>0$ there exists $\lambda_0>0$ such that for  any $\lambda>\lambda_0>0$ and $\vk{r}\in V^+$ and $\tz\in\mathcal{M}$
\begin{align}\label{eq}\left|\mathcal{H}_W^{p_{2,\vk{r}}^+g_{2,\vk{r}}^+(\bt_2)}- \mathcal{H}_W^{p_{2,\vk{r}}^+g_{2,\vk{r}}^+(\bt_2)}([0,\lambda]^{k_2})\lambda^{-k_1}\right|<\epsilon,
\left|\mathcal{H}_W^{p_2(\tz)g_2(\vk{a}_2^{-1}(\tz)\bt_2)}-\mathcal{H}_W^{p_2(\tz)g_2(\vk{a}_2^{-1}(\tz)\bt_2)}
([0,\lambda]^{k_2})\lambda^{-k_1}\right|<\epsilon.\end{align}

Hence it follows that, as $u\rw\IF$ and $\lambda>\lambda_0$,
\BQNY
\frac{\sum_{\vk{r}\in V^-}\pu{E_{1,\vk{r}}(u)}}{\Theta_1(u)}&\geq& \sum_{\vk{r}\in V^-}\mathcal{H}_W^{p_{2,\vk{r}}^+g_{2,\vk{r}}^+(\bt_2)}\prod_{i=1}^{k_1} a_{i,\vk{r}}^-\int_{\bt_1\in[0,\IF)^{k_1-k}} e^{-p_{1,\vk{r}}^+ g_1(\bt)}d\bt_1v^{k_0}\nonumber\\
&\geq&\int_{ \mathcal{M}}\sum_{\vk{r}\in  V^-}\left((\mathcal{H}_W^{p_{2,\vk{r}}^+g_{2,\vk{r}}^+(\bt_2)}([0,\lambda]^{k_1})\lambda^{-k_1}-\epsilon)\prod_{i=1}^{k_1} a_{i,\vk{r}}^-\int_{\bt_1\in[0,\IF)^{k_1-k}} e^{-p_{1,\vk{r}}^+ g_1(\bt)}d\bt_1\right) \mathbb{I}_{\mathcal{M}_{\vk{r}}}(\tz)d\tz.
\EQNY
Note that for any fixed $\tz\in\mathcal{M}^o$, where $\mathcal{M}^o\subset \mathcal{M}$ is the interior of $\mathcal{M}$,
\BQNY
&&\lim_{v\rw 0}\sum_{\vk{r}\in V^-}\left((\mathcal{H}_W^{p_{2,\vk{r}}^+g_{2,\vk{r}}^+(\bt_2)}([0,\lambda]^{k_1})\lambda^{-k_1}-\epsilon)\prod_{i=1}^{k_1} a_{i,\vk{r}}^-\int_{\bt_1\in[0,\IF)^{k_1-k}} e^{-p_{1,\vk{r}}^+ g_1(\bt)}d\bt_1\right)\mathbb{I}_{\mathcal{M}_{\vk{r}}}(\tz)\\
&&\quad = (\mathcal{H}_W^{p_2(\tz)g_2(\vk{a}_2^{-1}(\tz)\bt_2)}([0,\lambda]^{k_1})\lambda^{-k_1}-\epsilon) \LT(\prod_{i=1}^{k_1}a_i(\tz)\RT) \int_{\bt_1\in[0,\IF)^{k_1-k}} e^{-p_1(\tz)g_1(\bt_1)}d\bt_1\\
&&\quad \geq  (\mathcal{H}_W^{p_2(\tz)g_2(\vk{a}_2^{-1}(\tz)\bt_2)}-2\epsilon) \LT(\prod_{i=1}^{k_1}a_i(\tz)\RT) \int_{\bt_1\in[0,\IF)^{k_1-k}} e^{-p_1(\tz)g_1(\bt_1)}d\bt_1\\
&&\quad \geq \mathcal{H}_W^{p_2(\tz)g_2(\vk{a}_2^{-1}(\tz)\bt_2)} \LT(\prod_{i=1}^{k_1}a_i(\tz)\RT) \int_{\bt_1\in[0,\IF)^{k_1-k}} e^{-p_1(\tz)g_1(\bt_1)}d\bt_1,\quad \epsilon \rw 0.
\EQNY
Moreover, it is clear that there exits $\mathbb{Q}<\IF$ such that for any $\lambda>1$ and $v>0$
$$\left((\mathcal{H}_W^{p_{2,\vk{r}}^+g_{2,\vk{r}}^+(\bt_2)}([0,\lambda]^{k_1})\lambda^{-k_1}-\epsilon)\prod_{i=1}^{k_1} a_{i,\vk{r}}^-\int_{\bt_1\in[0,\IF)^{k_1-k}} e^{-p_{1,\vk{r}}^+ g_1(\bt)}d\bt_1\right) \mathbb{I}_{\mathcal{M}_{\vk{r}}}<\mathbb{Q}_8.$$
Consequently, dominated convergence theorem gives
\BQN\label{eq1}
\liminf_{u\rw\IF}\frac{\sum_{\vk{r}\in V^-}\pu{E_{1,\vk{r}}(u)}}{\Theta_1(u)}
\geq \int_{\mathcal{M}}
\left(\mathcal{H}_W^{p_2(\tz)g_2(\vk{a}_2^{-1}(\tz)\bt_2)} \LT(\prod_{i=1}^{k_1}a_i(\tz)\RT) \int_{\bt_1\in[0,\IF)^{k_1-k}} e^{-p_1(\tz)g_1(\bt_1)}d\bt_1\right)d\tz.\nonumber\\
\EQN
\cLa{Next we focus on  the double-sum term in (\ref{double}). For $\vk{r}\in V^-, \vk{r}'\in V^-, M_{\vk{r}}\cap M_{\vk{r}'}=\emptyset$,} we have
\begin{align*}
\pu{E_{1,\vk{r}}(u), E_{1,\vk{r}'}(u)}\leq \mathbb{P}\left(\sup_{\s\in E_{1,\vk{r}}, \t\in E_{1,\vk{r}'}} X(\s)+X(\t)>2u\right).
\end{align*}
By {\bf A1} and (\ref{reql}), there exists $0<\delta<1$ such that for all $\vk{r}\in V^-, \vk{r}'\in V^-, M_{\vk{r}}\cap M_{\vk{r}'}=\emptyset$,
\begin{align*}
\sup_{\s\in E_{1,\vk{r}}, \t\in E_{1,\vk{r}'}} Var(X(\s)+X(\t))<4-\delta.
\end{align*}
By Borell-TIS inequality (see, e.g., Theorem 2.1.1 in \cite{AdlerTaylor}), we have for $u>a$
\begin{align*}
\mathbb{P}\left(\sup_{\s\in E_{1,\vk{r}}, \t\in E_{1,\vk{r}'}} X(\s)+X(\t)>2u\right)\leq e^{-\frac{4(u-a)^2}{2(4-\delta)}},
\end{align*}
where $a=\frac{\mathbb{E}\left(\sup_{\s\in \mathcal{A}, \t\in \mathcal{A}} X(\s)+X(\t)\right)}{2}=\mathbb{E}(\sup_{\t\in \mathcal{A}}X(\t))$.
Consequently,
\begin{align*}
\sum_{\vk{r}, \vk{r}'\in V^-, M_{\vk{r}}\cap M_{\vk{r}'}= \emptyset}\pu{E_{1,\vk{r}}(u), E_{1,\vk{r}'}(u)}\leq \mathbb{Q}  e^{-\frac{4(u-a)^2}{2(4-\delta)}}=o(\Theta_1(u)),~~u\to\infty.
\end{align*}
For $\vk{r}, \vk{r}'\in V^-, \vk{r}\neq \vk{r}', M_{\vk{r}}\cap M_{\vk{r}'}\neq\emptyset$,
\begin{align*}
\pu{E_{1,\vk{r}}(u), E_{1,\vk{r}'}(u)}=\pu{E_{1,\vk{r}}(u)}+\pu{E_{1,\vk{r}'}}-\pu{E_{1,\vk{r}}(u), E_{1,\vk{r}'}(u)}.
\end{align*}
Hence in light of arguments of (\ref{eq}) and (\ref{eq1}), we have
\begin{align*}
\sum_{\vk{r}, \vk{r}'\in V^-, \vk{r}\neq \vk{r}', M_{\vk{r}}\cap M_{\vk{r}'}\neq\emptyset}\pu{E_{1,\vk{r}}(u), E_{1,\vk{r}'}(u)}=o(\Theta_1(u)),~u\to\infty, v\to 0.
\end{align*}
Therefore we have
\begin{align*}
\sum_{\vk{r}, \vk{r}'\in V^-, \vk{r}\neq \vk{r}'}\pu{E_{1,\vk{r}}(u), E_{1,\vk{r}'}(u)}=o(\Theta_1(u)),~u\to\infty, v\to 0,
\end{align*}
implying
\begin{align*}
\liminf_{u\rw\IF}\frac{\pu{E_{1}(u)}}{\Theta_1(u)}
\geq \int_{\mathcal{M}}
\left(\mathcal{H}_W^{p_2(\tz)g_2(\vk{a}_2^{-1}(\tz)\bt_2)} \LT(\prod_{i=1}^{k_1}a_i(\tz)\RT) \int_{\bt_1\in[0,\IF)^{k_1-k}} e^{-p_1(\tz)g_1(\bt_1)}d\bt_1\right)d\tz.
\end{align*}
Similarly,
\BQNY
\limsup_{u\rw\IF}\frac{\pu{E_2(u)}}{\Theta_1(u)}\leq\int_{\mathcal{M}}\left(\mathcal{H}_W^{p_2(\tz)g_2(\vk{a}_2^{-1}(\tz)\bt_2)} \LT(\prod_{i=1}^{k_1}a_i(\tz)\RT) \int_{\bt_1\in[0,\IF)^{k_1-k}} e^{-p_1(\tz)g_1(\bt_1)}d\bt_1\right) d\tz , v\rw 0.
\EQNY
Inserting the above two inequalities and (\ref{err0}) in (\ref{boundlem1}), we establish the claim.
\QED

\section{Proof of Proposition \ref{ThmM1}}\label{s.Prop31}
We distinguish three scenarios according to the value of $\alpha$: $\alpha\in (0,1)$, $\alpha=1$ and $\alpha\in (1,2)$.
The cases of $\alpha\in(0,1)$ and $\alpha\in(1,2)$  \cLa{can be derived from} \cite{Pit96}[Theorem 8.2],
where the maximum of the variance is attained at finite number of points.
However, the case $\alpha=1$ is essentially  different from the aforementioned two cases
in the sense that depending on $a_i$ the maximum of the variance can be achieved at a
set of positive Lebesgue measure of dimension
$\mathfrak{m}-1$, where $\mathfrak{m}$ is defined in (\ref{mmd}).
 We apply Theorem \ref{MThm2} to this case.

For $Z^{\alpha}(\t)$ \cLa{introduced} in \eqref{ZZ1} with $\alpha\in(0,2)$, we write $\sigma_Z^2$
for the variance of $Z^{\alpha}$ and $r_Z$ for its correlation function.
Moreover, we denote $\sigma_*= \max_{\t\in \Sn} \sigma_Z(\vk t)$ and recall
that $\mathcal{S}_{n}=\{0=t_0\leq t_1\leq\cdots\leq t_n\leq t_{n+1}=1\}$.
The expansions of  $\sigma_Z$ and $r_Z$ are displayed in the following lemma
which \cLa{is} crucial for the proof of Proposition \ref{ThmM1}.
\kk{We skip  its  proof  as it only needs some standard but tedious calculations.}

\BL\label{lem1}
i) For $\alpha \in (0,1)$,  the standard deviation $\sigma_Z$  attains its maximum  on $\mathcal{S}_{n}$ at only one point
$\vk z_0= (z_1\ldot z_{n})\in \mathcal{S}_{n}$ with
$
z_i=\frac{\sum_{j=1}^ia_j^{\frac{2}{1-\alpha}}}
{\sum_{j=1}^{n+1}a_j^{\frac{2}{1-\alpha}}},\ i=1\ldot n,
$
and its maximum value is
$
\sigma_*=\LT(\sum_{i=1}^{n+1}a_i^{\frac{2}{1-\alpha}}\RT)^{\frac{1-\alpha}{2}}.
$
Moreover,
\BQN\label{var1}
\lim_{\delta\rw0}\underset{\abs{\t-\z_0}\leq\delta}{\sup_{\t\in\Sn}} \abs{\frac{1-\frac{\sigma_Z(\t)}{\sigma_*}}{\frac{\alpha(1-\alpha)\LT(\sum_{i=1}^{n+1}a_i^{\frac{2}{1-\alpha}}\RT)}{4}
		\sum_{i=1}^{n+1}a_i^{\frac{2}{\alpha-1}}
		\LT((t_i-z_i)-(t_{i-1}-z_{i-1})\RT)^2}-1}=0
\EQN
with \cLa{$z_0:=0, z_{n+1}:=1$,} and
\BQN\label{r1}
\lim_{\delta\rw 0}\underset{\abs{\s-\z_0},\abs{\t-\z_0} <\delta }{\sup_{\s\neq \t, \s, \t\in\Sn}} \abs{\frac{1-r_Z(\s,\t)}{\frac{1}{2\sigma_*^2}\LT(\sum_{i=1}^{n}(a_i^2+a_{i+1}^2)
		\abs{s_i-t_i}^\alpha\RT)}-1}=0.
\EQN

ii) For $\alpha=1$ \cLa{and $\mathfrak{m}$ defined in \eqref{mmd}}, if $\mathfrak{m}=n+1$, $\sigma_Z(\t)\equiv1,\ \t \in \Sn$, and if $\mathfrak{m}<n+1$,  function $\sigma_Z$  attains its maximum equal to $1$ on $\mathcal{S}_{n}$ at
$\mathcal{M}=\{\t\in\Sn:\sum_{j\in \mathcal{N}}\abs{t_j-t_{j-1}}=1\}$ and satisfies
\BQN\label{sigma21}
\lim_{\delta\rw0}\sup_{\z\in \mathcal{M}}\underset{\abs{\t-\z}\leq\delta}{\sup_{\t\in\Sn}}\abs{\frac{1-\sigma_Z(\t)}
{\frac{1}{2}\sum_{j\in\MNB}(1-a_j^2)\abs{t_j-t_{j-1}}}-1}=0.
\EQN
In addition, for $1\leq \mathfrak{m}\leq n+1$, we have
\BQN\label{r2}
\lim_{\delta\rw 0}\sup_{\z\in \mathcal{M}}\underset{\abs{\s-\z}, \abs{\t-\z}<\delta }{\sup_{\s\neq \t, \s, \t\in\Sn}}\left|\frac{1-r_Z(\s,\t)}{\frac{1}{2}\sum_{i=1}^{n+1}
a_i^2\min\LT(\abs{t_{i-1}-s_{i-1}}+\abs{t_i-s_i},\abs{t_i-t_{i-1}}+\abs{s_i-s_{i-1}}\RT)}-1\right|=0.
\EQN

iii) For $\alpha\in(1,2)$,  function $\sigma_Z $  attains it maximum on $\Sn$ at $\mathfrak{m}$ points $\z^{(j)},\ j\in\mathcal{N}=\{i:a_i=1, i=1\ldot n+1\}$,
where $\z^{(j)}=(0\ldot 0, 1, 1\ldot 1)$ (the first $1$ stands at the $j$-th coordinate) if $j\in\mathcal{N}$ and $j<n+1$, and $\z^{(n+1)}=(0\ldot 0)$ if $n+1\in\mathcal{N}$. We further  have that $\sigma_*=1$ and  as $\t \to \vk z^{(j)}$
\BQN\label{var3}
\lim_{\delta\rw0}\underset{\abs{\t-\z^{(j)}}\leq\delta}{\sup_{\t\in\Sn}}
\abs{\frac{1-\sigma_Z(\t)}{\frac{1}{2}\LT(\alpha \abs{t_j-t_{j-1}-1}-
		\sum_{1 \le i \le n+1, i\not= j }a_i^2\abs{ t_i-t_{i-1}}^\alpha\RT)}-1}=0.
\EQN
\EL

{\bf Case 1. $\alpha\in(0,1)$:}   From \nelem{lem1} i), we have \cLa{that} $\sigma_Z$  on $\mathcal{S}_{n}$ attains its maximum $\sigma_*$ at the unique point
$\vk z_0= (z_1\ldot z_{n})$ with
\BQNY
z_i=\frac{\sum_{j=1}^ia_j^{\frac{2}{1-\alpha}}}
{\sum_{j=1}^{n+1}a_j^{\frac{2}{1-\alpha}}},\ i=1\ldot n.
\EQNY
Moreover, from \eqref{var1} we have  for $\t\in\Sn$
\COM{\BQNY
1-\frac{\sigma_Z(\t)}{\sigma_*}=\frac{\alpha(1-\alpha)\LT(\sum_{i=1}^{n+1}a_i^{\frac{2}{1-\alpha}}\RT)}{4}
		\sum_{i=1}^{n+1}a_i^{\frac{2}{\alpha-1}}
		\LT((t_i-z_i)-(t_{i-1}-z_{i-1})\RT)^2(1+o(1)),~\t\to\z_0.
\EQNY
Is it better to write it in the following way since there in no definition for $z_0$ and $z_{n+1}$?}
\cLa{\BQNY
1-\frac{\sigma_Z(\t)}{\sigma_*}\sim\frac{\alpha(1-\alpha)\LT(\sum_{i=1}^{n+1}a_i^{\frac{2}{1-\alpha}}\RT)}{4}
		\left(a_1^{\frac{2}{\alpha-1}}
		(t_1-z_1)^2+a_{n+1}^{\frac{2}{\alpha-1}}
		(t_n-z_n)^2+\sum_{i=2}^{n}a_i^{\frac{2}{\alpha-1}}
		\LT((t_i-z_i)-(t_{i-1}-z_{i-1})\RT)^2\right),
\EQNY}

as $\abs{\t-\z_0}\rw 0$ and from \eqref{r1} for $\t, \s\in\Sn$
\BQNY
1-r_Z(\s,\t)\sim\frac{1}{2\sigma_*^2}\LT(\sum_{i=1}^{n}(a_i^2+a_{i+1}^2)
		\abs{s_i-t_i}^\alpha\RT),
\EQNY
as $\abs{\s-\z_0},\abs{\t-\z_0} \rw 0$.
Further, we have
\BQN\label{eq10}
\E{(Z^{\alpha}(\s)-Z^{\alpha}(\t))^2}
\leq 4\sum_{i=1}^{n}
\abs{t_i-s_i}^\alpha.
\EQN
Thus by \cite{Pit96} [Theorem 8.2] we obtain that, as $u\rw\IF$,
\BQNY
\pk{\sup_{\t\in\mathcal{S}_{n}}Z^\alpha(\t)>u}\sim \LT( \mathcal{H}_{B^{\alpha}}\RT)^n\prod_{i=1}^n\LT(\frac{a_i^2+a_{i+1}^2}{2\sigma^2_*}\RT)^{1/\alpha}
	\LT(\frac{u}{\sigma_*}\RT)^{(2/\alpha-1)n}\int_{\R^n}e^{-f(\x)}d\x\Psi\LT(\frac{u}{\sigma_*}\RT),
\EQNY
where
\BQNY
f(\x)=\frac{\alpha(1-\alpha)\LT(\sum_{i=1}^{n+1}a_i^{\frac{2}{1-\alpha}}\RT)}{4}
		\left(a_1^{\frac{2}{\alpha-1}}
		x_1^2+a_{n+1}^{\frac{2}{\alpha-1}}
		x_n^2+\sum_{i=2}^{n}a_i^{\frac{2}{\alpha-1}}
		\LT(x_i-x_{i-1}\RT)^2\right),~\x\in\mathbb{R}^n.
\EQNY
Direct calculation shows
\BQNY
\int_{\R^n}e^{-f(\x)}d\x=\LT(\frac{4\pi}{\alpha(1-\alpha)}\RT)^{\frac{n}{2}}\sigma_*^{-\frac{n}{1-\alpha}}\left(\sum_{j=1}^{n+1}\prod_{i\neq j}a_i^{\frac{2}{\alpha-1}}\right)^{-\frac{1}{2}}.
\EQNY
This completes the proof of this case.

{\bf Case 2. $\alpha=1$:}
First we consider the case $\mathfrak{m}<n+1$. Let  $k^*=\max\{i\in\mathcal{N}\}$ and denote
\BQNY
\mathcal{N}_0=\{i\in\mathcal{N}, i<k^*\},\
\MNB_0=\{i\in\MNB, i<k^*\}.
\EQNY
In order to facilitate our analysis, we make the following transformation:
\BQNY
x_i=t_i,\ i\in\mathcal{N}_0,\
x_i=t_i-t_{i-1},\ i\in\MNB,
\EQNY
implying that $\x=(x_1\ldot x_{k^*-1},x_{k^*+1}\ldot x_{n+1})\in[0,1]^n$ and
\BQN\label{tt}
t_i=t_i(\x)=\LT\{
\begin{array}{ll}
x_i,& \text{if}\ i\in \mathcal{N}_0,\\
1-\sum_{j=i+1}^{n+1}x_j,& \text{if}\ i\geq k^*,\\
\sum_{j=\max\{k\in\mathcal{N}:k<i\}}^{i}x_j,& \text{if}\ i\in\MNB_0,
\end{array}
\RT.
\EQN
with the convention that $\max\emptyset=0$.
Define $Y(\x)=Z(\t(\x))$ and $\widetilde{\mathcal{S}}_n=\{\x:\t(\x)\in\mathcal{S}_n\}$ with $\t(\x)$ given in \eqref{tt}. By \nelem{lem1} ii) it follows that
$\sigma_{Y}(\x)$, the standard deviation of $Y(\x)$, attains its maximum equal to $1$ at
$$\{\x\in\widetilde{\mathcal{S}}_n: x_i=0,\ \text{if}\ i\in\MNB\} .$$

Moreover,  let $\widetilde{\x}=(x_i)_{i\in \mathcal{N}_0 }$, $\bar{\x}=(x_i)_{i\in\mathcal{N}^c}$  and denote for any $\delta\in (0, \frac{1}{(n+1)^2})$
\BQNY\label{eq:M}
\widetilde{\mathcal{S}}^*_n(\delta)&=&\LT\{\x\in\widetilde{\mathcal{S}}_n:0\leq x_i\leq \frac{\delta}{(n+1)^2}, \text{if} \ i\in\MNB\RT\},\nonumber\\
\widetilde{\mathcal{M}}&=&\{\widetilde{\x}\in[0,1]^{\mathfrak{m}-1}:  x_i\leq x_j,\quad  \text{if}\ i,j\in \mathcal{N}_0\ \text{and}\ i<j\},\nonumber\\
\widetilde{\mathcal{M}}(\delta)&=&\{\widetilde{\x}\in[\delta,1-\delta]^{\mathfrak{m}-1}:  x_j-x_i\geq \delta,\ \text{if}\ i,j\in \mathcal{N}_0\ \text{and}\ i<j\}\subseteq\widetilde{\mathcal{M}},\\
\widetilde{\mathcal{S}}_n(\delta)&=&\{\x\in\widetilde{\mathcal{S}}^*_n(\delta): \widetilde{\x}\in\widetilde{\mathcal{M}}(\delta)\}.\nonumber
\EQNY
We notice  that
\BQN
\pk{\sup_{\x\in \widetilde{\mathcal{S}}_n}Y(\x)>u}\geq \pk{\sup_{\x\in \widetilde{\mathcal{S}}_n(\delta)}Y(\x)>u},\label{bounduplow1}
\EQN
and
\BQN\label{ineq}
\pk{\sup_{\x\in \widetilde{\mathcal{S}}_n}Y(\x)>u}&\leq&
\pk{\sup_{\x\in \widetilde{\mathcal{S}}_n\setminus\widetilde{\mathcal{S}}^*_n(\delta)}Y(\x)>u}\nonumber\\
&&+\pk{\sup_{\x\in\widetilde{\mathcal{S}}^*_n(\delta)\setminus\widetilde{\mathcal{S}}_n(\delta)}Y(\x)>u}+\pk{\sup_{\x\in \widetilde{\mathcal{S}}_n(\delta)}Y(\x)>u}.\label{boundup1}
\EQN
The idea  of the proof is first to apply Theorem \ref{MThm2} to obtain the asymptotics of $\pk{\sup_{\x\in\widetilde{\mathcal{S}}_n(\delta)}Y(\x)>u}$ as $u\to\infty$ and then to show that \cLa{the} other two terms in (\ref{ineq}) are asymptotically negligible.  Let us begin with finding the asymptotics of $\pk{\sup_{\x\in\widetilde{\mathcal{S}}_n(\delta)}Y(\x)>u}$. First observe $$\widetilde{\mathcal{S}}_n(\delta)=\{\x: \widetilde{\x}\in\widetilde{\mathcal{M}}(\delta),~0\leq x_i\leq \frac{\delta}{(n+1)^2}, \text{if} \ i\in\MNB\},$$ which is a set satisfying the assumption in Theorem \ref{MThm2}. Moreover,
 It follows from \eqref{sigma21} that
\BQN\label{eqcona1}
\lim_{\delta\rw 0}\sup_{\x\in \widetilde{\mathcal{S}}_n(\delta)}\abs{\frac{1-\sigma_{Y}(\x)}{\frac{1}{2}\sum_{i\in\MNB}(1-a_i^2)x_i}-1}=0.
\EQN
\cLa{Taking $\tt=\widetilde{\x}$  and $\bt_2=\bar{\x}$} in \netheo{MThm2}, \eqref{eqcona1} implies that {\bf A2} holds with $g_2(\bar{\x})=\frac{1}{2}\sum_{i\in\MNB}(1-a_i^2)x_i$  and $p_2(\widetilde{\x})=1$ for $\widetilde{\x}\in\widetilde{\mathcal{S}}^*_n(\delta)$.
\cLa{We note that  $\Lambda_1=\Lambda_3=\emptyset$ in this case.}

We next check Assumption {\bf A1}. To compute the correlation, we need the following inequalities.  Note that
for $\x,\y\in \widetilde{\mathcal{S}}_n(\delta)$ and  $\abs{\x-\y}<\frac{\delta}{(n+1)^2}$, if $i\in\mathcal{N}_0$,
$$\abs{x_i-y_i}+\abs{t_{i-1}(\x)-t_{i-1}(\y)}< \frac{\delta}{(n+1)^2}+\frac{n\delta}{(n+1)^2}=\frac{\delta}{n+1}\leq\frac{\delta}{2}$$
and
\BQNY
\abs{t_i(\x)-t_{i-1}(\x)}=\LT\{
\begin{array}{ll}
 \abs{x_i-x_{i-1}}\geq\delta & \text{if}\  i-1\in\mathcal{N}_0, \\
 \abs{x_i-\sum_{j=\max\{k\in\mathcal{N}:k<i-1\}}^{i-1}x_j}\geq\delta-\frac{n\delta}{(n+1)^2}> \frac{\delta}{2}& \text{if}\  i-1\in\mathcal{N}^c,
\end{array}
\RT.
\EQNY
if $i=k^*$,
$$\abs{t_{k^*-1}(\y)-t_{k^*-1}(\x)}+\abs{t_{k^*}(\y)-t_{k^*}(\x)}<\frac{n\delta}{(n+1)^2}<\frac{\delta}{2},$$
and
\BQNY
\abs{t_{k^*}(\x)-t_{k^*-1}(\x)}=\LT\{
\begin{array}{ll}
 \abs{1-\sum_{j=k^*+1}^{n+1}x_j-x_{k^*-1}}\geq 1-(1-\delta)-\frac{n\delta}{(n+1)^2}> \frac{\delta}{2} & \text{if}\  k^*-1\in\mathcal{N}_0, \\
 \abs{1-\sum_{j=k^*+1}^{n+1}x_j-\sum_{j=\max\{k\in\mathcal{N}:k<k^*-1\}}^{k^*-1}x_j}\geq 1-(1-\delta)-\frac{n\delta}{(n+1)^2}> \frac{\delta}{2}& \text{if}\  k^*-1\in\mathcal{N}^c.
\end{array}
\RT.
\EQNY
Hence for $r_Y(\x,\y)$, the correlation function of $Y(\x)$, we derive from  ii) of \nelem{lem1} that  for $\x,\y\in \widetilde{\mathcal{S}}_n(\delta)$ and $\abs{\x-\y}<\frac{\delta}{(n+1)^2}$, as $\delta\rw 0$
\cLa{\begin{align}
1-r_Y(\x,\y)=&1-r_Z(\t(\x),\t(\y))\nonumber\\
\sim&\frac{1}{2}\sum_{i=1}^{n+1}
a_i^2\min\LT(\abs{t_{i-1}(\y)-t_{i-1}(\x)}+\abs{t_i(\y)-t_i(\x)},\abs{t_i(\y)-t_{i-1}(\y)}
+\abs{t_i(\x)-t_{i-1}(\x)}\RT)\nonumber\\
=&\frac{1}{2}\sum_{i\in\mathcal{N}}(\abs{t_{i-1}(\y)-t_{i-1}(\x)}
+\abs{t_i(\y)-t_{i}(\x)})\nonumber\\
&+\frac{1}{2}\sum_{i\in \mathcal{N}^c}
a_i^2\min\LT(\abs{t_{i-1}(\y)-t_{i-1}(\x)}+\abs{t_i(\y)-t_i(\x)},\abs{t_i(\y)-t_{i-1}(\y)}
+\abs{t_i(\x)-t_{i-1}(\x)}\RT)\nonumber\\
=&\frac{1}{2}\sum_{i\in\mathcal{N}_0}
\LT(\abs{x_i-y_i}+\abs{t_{i-1}(\x)-t_{i-1}(\y)}\RT)\nonumber\\
&+
\frac{1}{2}\abs{t_{k^*-1}(\x)-t_{k^*-1}(\y)}
+\frac{1}{2}\abs{\sum_{j=k^*+1}^{n+1}(x_j-y_j)}\nonumber\\
&+\frac{1}{2}\sum_{i\in\MNB_0}
a_i^2\min\LT(\abs{t_{i-1}(\x)-t_{i-1}(\y)}+\abs{t_i(\x)-t_i(\y)},x_i+y_i\RT)\nonumber\\
&+\frac{1}{2}\sum_{i=k^*+1}^{n+1}
a_i^2\min\LT(\abs{\sum_{j=i}^{n+1}(x_j-y_j)}+\abs{\sum_{j=i+1}^{n+1}(x_j-y_j)},x_i+y_i\RT). \label{yrr}
\end{align}}

By \eqref{tt}, we have for any $i=1,\cdots,n+1$
\BQNY
\abs{t_i(\y)-t_i(\x)}\leq \underset{ i\neq k^*}{\sum_{i=1}^{n+1}}|x_i-y_i|.
\EQNY
Then for $\x,\y\in \widetilde{\mathcal{S}}_n(\delta)$  and $|\x-\y|<\frac{\delta}{(n+1)^2}$ with $\delta>0$ sufficiently small
\BQNY
\frac{1}{2}\sum_{i\in \mathcal{N}_0}|x_i-y_i|\leq 1-r_Y(\x,\y)\leq \mathbb{Q}\underset{ i\neq k^*}{\sum_{i=1}^{n+1}}|x_i-y_i|,
\EQNY
implying that (\ref{boundrr1}) holds.

Recall
\begin{align}\label{eqW}
W(\x)&=\frac{\sqrt{2}}{2}\sum_{i\in\mathcal{N}}
\LT(B_i(s_i(\x))-\widetilde{B}_i(s_{i-1}(\x))\RT)+\frac{\sqrt{2}}{2}\sum_{i\in\MNB}
a_i\LT(B_{i}(s_i(\x))-B_i(s_{i-1}(\x))\RT),
\end{align}
where $B_i, \widetilde{B}_i$ are iid standard Brownian motions and
\BQNY
s_i(\x)=\LT\{
\begin{array}{ll}
x_i,& \text{if}\ i\in \mathcal{N}_0,\\
\sum_{j=\max\{k\in\mathcal{N}:k<i\}}^{i}x_j,& \text{if}\ i\in\MNB_0,\\
\sum_{j=i+1}^{n+1}x_j,& \text{if}\ i\geq k^*,
\end{array}
\RT.
\EQNY with the convention that $\max\emptyset=0$.
Direct calculation gives that  $\E{\LT(W(\x)-W(\y)\RT)^2}$ coincides with \eqref{yrr} for any $\x,\y\in[0,\IF)^n$.\\
This implies that (\ref{Rr}) holds with $W$ given in (\ref{eqW}) and $\vk{a}(\tilde{x})\equiv 1$  for $\tilde{x}\in\widetilde{\mathcal{M}}(\delta)$.

Using \eqref{yrr} and the fact that for any $i=1,\cdots,n$, $s_i(\x)-s_i(\y)$ is the absolute value of the combination of $x_j-y_j, \ j\in\{1,\cdots,k^*-1,k^*+1,\cdots, n+1\}$,  we derive that for a fixed $\bar{\x}$ the increments of $W(\x)=W(\widetilde{\x},\bar{\x})$ are homogeneous with respect to $\widetilde{\x}$. In addition, it is easy to check that (\ref{reql}) also holds.  Hence {\bf A1} is satisfied.

Consequently, by \netheo{MThm2},  as $u\rw\IF$, we have
\BQN\label{asyp1}
\pk{\sup_{\x\in\widetilde{\mathcal{S}}_n(\delta)}Y(\x)>u}\sim v_{\mathfrak{m}-1}\LT(\widetilde{\mathcal{M}}(\delta)\RT)\mathcal{H}_Wu^{2(\mathfrak{m}-1)}\Psi(u),
\EQN
where
\BQNY
\mathcal{H}_W&=&\lim_{\lambda\rw\IF}\frac{1}{\lambda^{\mathfrak{m}-1}}
\E{\sup_{\x\in[0,\lambda]^n}e^{ \sqrt{2}W(\x)-\sigma_W^2(\x)-\frac{1}{2}\sum_{j\in\MNB}(1-a_j^2)x_j}}\\
&=&\lim_{\lambda\rw\IF}\frac{1}{\lambda^{\mathfrak{m}-1}}
\E{\sup_{\x\in[0,\lambda]^n}e^{ \sqrt{2}W(\x)-\LT(\underset{i\neq k^*}{\sum_{i=1}^{n+1}}x_i\RT)}}.
\EQNY

 We now proceed to the negligibility of the other two terms in (\ref{boundup1}).
In light of  Borell-TIS inequality\COM{(see, e.g., Theorem 2.1.1 in \cite{AdlerTaylor}) }, we have
\BQN\label{neg22}
\pk{\sup_{\x\in   \widetilde{\mathcal{S}}_n\setminus\widetilde{\mathcal{S}}^*_n(\delta) }Y(\x)>u}
\leq \exp\left(\frac{(u-\mathbb{E}(\sup_{\x\in   \widetilde{\mathcal{S}}_n\setminus\widetilde{\mathcal{S}}^*_n(\delta) }Y(\x)))^2}{2(1-\epsilon)^2}\right)=o(\Psi(u)),~u\to\infty,
\EQN
where $\vn=1-\sup_{\x\in \widetilde{\mathcal{S}}_n\setminus\widetilde{\mathcal{S}}^*_n(\delta)}\sigma_Y(\x).$
By Slepian's inequality and Theorem \ref{MThm2}, we have
\BQN\label{neg1}
\pk{\sup_{\x\in\widetilde{\mathcal{S}}^*_n(\delta)\setminus\widetilde{\mathcal{S}}_n(\delta)}Y(\x)>u}
&\leq& v_{\mathfrak{m}-1}\LT(\widetilde{\mathcal{M}}\setminus\widetilde{\mathcal{M}}(\delta)\RT)
\widetilde{\mathcal{H}}_{W_1}u^{2(\mathfrak{m}-1)}\Psi(u)\nonumber\\
&=&o\LT(u^{2(\mathfrak{m}-1)}\Psi(u)\RT),\ u\rw\IF,\  \delta\rw 0.
\EQN
Combination of the fact that
\BQNY
\lim_{\delta\to 0}v_{\mathfrak{m}-1}\LT(\widetilde{\mathcal{M}}(\delta)\RT)=v_{\mathfrak{m}-1}\LT(\widetilde{\mathcal{M}}\RT)=\frac{1}{(\mathfrak{m}-1)!}
\EQNY
with 
 \eqref{bounduplow1}, \eqref{boundup1}, and (\ref{asyp1})-\eqref{neg1} leads to
\BQNY
\pk{\sup_{\t\in\mathcal{S}_{n}}Z(\t)>u}=\pk{\sup_{\x\in\widetilde{\mathcal{S}}_n}Y(\x)>u}\sim \frac{1}{(\mathfrak{m}-1)!}\mathcal{H}_Wu^{2(\mathfrak{m}-1)}\Psi(u),\ u\rw\IF.
\EQNY

\underline{Case $\mathfrak{m}=n+1$}:
For some small $\vn\in(0,1)$, define $E(\vn)=\{\t\in \Sn: t_i-t_{i-1}\geq \vn, i=1\ldot n+1\}$. Then we have
\BQN\label{boundul}
\pk{\sup_{\t\in E(\vn)}Z(\t)>u}\leq\pk{\sup_{\t\in \Sn}Z(\t)>u}\leq \pk{\sup_{\t\in \Sn\setminus E(\vn)}Z(\t)>u}+\pk{\sup_{\t\in E(\vn)}Z(\t)>u}.
\EQN
Let us first derive the asymptotics of  $Z$ over $E(\vn)$.
For $\s, \t\in E(\vn)$, by \eqref{r2} we have
\BQNY
1-r(\s,\t)\sim\sum_{i=1}^{n}\abs{s_i-t_i}, \abs{\t-\s}\rw 0.
\EQNY
Moreover, it follows \cLa{straightforwardly} that $Var(Z(\t))=1$ for $\t\in E(\vn)$ and $Corr(Z(\t), Z(\s))<1$ for any $\s\neq \t$ and $\s,\t\in E(\vn)$. Hence
by \cite{Pit96} [Lemma 7.1], we have
\BQN\label{2lb}
\pk{\sup_{\t\in E(\vn)}Z(\t)>u}\sim v_n(E(\vn))u^{2n}\Psi(u)\sim v_n(\Sn)u^{2n}\Psi(u) , \ u\rw\IF,\ \vn\rw 0.
\EQN

Moreover, by Slepian's inequality and \cite{Pit96} [Lemma 7.1],
\BQN\label{2err}
\pk{\sup_{\t\in \Sn\setminus E(\vn)}Z(\t)>u}
\leq v_n( \Sn\setminus E(\vn))(2\mathcal{H}_{B^1}\mathbb{Q}_4)^n u^{2n}\Psi(u)
=o\LT(u^{2n}\Psi(u)\RT),\ u\rw\IF, \vn\rw 0.
\EQN
Inserting \eqref{2lb} and \eqref{2err} into \eqref{boundul}, we obtain
\BQNY
\pk{\sup_{\t\in \Sn}Z(\t)>u}\sim \frac{1}{n!}u^{2n}\Psi(u),\ u\rw\IF.
\EQNY
The claim is established by ii) of Remark \ref{prop41}.

{\bf Case 3. $\alpha\in(1,2)$:}
For $\s,\t\in\Sn$, one can easily check that
\BQNY
r_Z(\s,\t)=\frac{\E{Z^\alpha(\t)Z^\alpha(\s)}}{\sigma_Z(\t)\sigma_Z(\s)}
=\frac{\sum_{i=1}^{n+1}a_i^2\E{(B_i^{\alpha}(t_i)-B_i^{\alpha}(t_{i-1}))
(B_i^{\alpha}(s_i)-B_i^{\alpha}(s_{i-1}))}}{\sigma_Z(\t)\sigma_Z(\s)}<1
\EQNY
if $\s\neq\t$.
In light of \nelem{lem1} iii), $\sigma_Z$ attains its maximum at $\mathfrak{m}$ distinct points $\z^{(j)},j\in\mathcal{N}$.
Consequently,  by \cite{Pit96} [Corollary 8.2], we have
\BQNY
\pk{\sup_{\t\in\mathcal{S}_{n}}Z^\alpha(\t)>u}
\sim\sum_{j\in \mathcal{N}}\pk{\sup_{\t\in\Pi_{\delta, j}}Z^\alpha(\t)>u},\ u\rw\IF,
\EQNY
where $\Pi_{\delta, j}=\LT\{\t\in\Sn: \abs{\t-\z^{(j)}}\leq \frac{1}{3}\RT\}.$\\
Define $E_j(u):=\{\t\in\Pi_{\delta, j}:1-\LT(\frac{\ln u}{u}\RT)^2\leq t_j-t_{j-1}\leq 1\}\ni \z^j$. Observe that
\BQNY
\pk{\sup_{\t\in E_j(u)}Z^\alpha(\t)>u}\leq \pk{\sup_{\t\in\Pi_{\delta, j}}Z^\alpha(\t)>u}\leq\pk{\sup_{\t\in E_j(u)}Z^\alpha(\t)>u}
+\pk{\sup_{\t\in\Pi_{\delta, j}\setminus E_j(u)}Z^\alpha(\t)>u}.
\EQNY
We first find the exact asymptotics of $\pk{\sup_{\t\in E_j(u)}Z^\alpha(\t)>u}$ as $u\to\infty$.
Clearly, for any $u\in\mathbb{R}$,
\BQNY
 \pk{\sup_{\t\in E_j(u)}Z^\alpha(\t)>u}\geq \pk{Z^\alpha(\z^j)>u}=\Psi(u).
\EQNY
Moreover, for  $\s,\t\in \Sn $, there exists a constant $c>0$ such that $\inf_{\t\in\Sn}\sigma_Z(\t)\geq \frac{1}{\sqrt{2c}}$. Hence  in light of (\ref{eq10}) we have
\BQN\label{eq11}
1-r_Z(\s,\t)
\leq  4c\sum_{i=1}^{n}\abs{t_i-s_i}^\alpha.
\EQN
Let $U_2(\t), \t\in \R^n$ be a centered homogeneous Gaussian field with continuous trajectories, unit variance and correlation function $r_{U_2}(\s,\t)$ satisfying
\BQNY
r_{U_2}(\s,\t)=1-\exp\LT(8c\sum_{i=1}^{n}\abs{t_i-s_i}^\alpha\RT).
\EQNY
Set $\widetilde{E}_j(u)=[0,\vn_1 u^{-2/\alpha}]^{j-1}\times[1-\vn_1 u^{-2/\alpha},1]^{n-j+1}$ for some constant $\vn_1\in(0,1)$. Then it follows that $E_j(u)\subset \widetilde{E}_j(u)$ for $u$ sufficiently large. By Slepian's inequality and \cite{Pit96} [Lemma 6.1]
\BQNY
\pk{\sup_{\t\in E_j(u)}Z^\alpha(\t)>u}\leq \pk{\sup_{\t\in \widetilde{E}_j(u)}U_2(\t)>u}
\sim \LT(\mathcal{H}_{B^{\alpha}}[0,(8c)^{1/\alpha}\vn_1]\RT)^{n}\Psi(u)\sim \Psi(u),
\EQNY
as $u\rw\IF, \vn_1\rw 0$, where
$$
\lim_{\lambda\rw0}\mathcal{H}_{B^{\alpha}}[0, \lambda]=\lim_{\lambda\rw0}\mathbb{E}\left\{\sup _{t \in[0, \lambda]} e^{\sqrt{2} B^{\alpha}(t)-t^{\alpha}}\right\}=1.
$$
Consequently,
\BQN\label{eq13}\pk{\sup_{\t\in E_j(u)}Z^\alpha(\t)>u}\sim \Psi(u),~u\to\infty.
\EQN

Note that for $\t\in\mathcal{S}_{n}$
$$\underset{i\neq j}{\sum_{i=1}^{n+1}}a_i^2\abs{t_i-t_{i-1}}^{\alpha}\leq \abs{t_j-t_{j-1}-1}.$$
Hence, by \eqref{var3}, for $u$ sufficiently large,
\cLa{\BQN\label{eq12}
\sup_{\t\in \Pi_{\delta, j}\setminus E_j(u) }\sigma_Z(\t)
&\leq&\sup_{\t\in \Pi_{\delta, j}\setminus E_j(u) }\LT( 1-\frac{(1-\vn)(\alpha-1)}{2}\abs{t_j-t_{j-1}-1}\RT)\nonumber\\
&\leq&1-\frac{(1-\vn)(\alpha-1)}{2}\LT(\frac{\ln u}{u}\RT)^2,
\EQN}
where $\vn\in(0,1)$ is a constant. In light of (\ref{eq11}) and (\ref{eq12}),
by  \cite{Pit96} [Theorem 8.1] we have,  for $u$ sufficiently large,
\BQNY
\pk{\sup_{\t\in\Pi_{\delta, j}\setminus E_j(u)}Z^\alpha(\t)>u}
\leq\mathbb{Q}_9u^{2n/\alpha}\Psi\LT(\frac{u}{
1-\frac{(1-\vn)(\alpha-1)}{2}\LT(\frac{\ln u}{u}\RT)^2}\RT)
=o\LT(\Psi\LT(u\RT)\RT),\ u\rw\IF,
\EQNY
which combined with (\ref{eq13}) leads to
\BQNY
\pk{\sup_{\t\in\Pi_{\delta, i}}Z^\alpha(\t)>u}\sim \pk{\sup_{\t\in E_j(u)}Z^\alpha(\t)>u}\sim \Psi(u),~u\to\infty.
\EQNY
Consequently, with $\mathfrak{m}=\sharp\mathcal{N}$, we obtain
\BQNY
\pk{\sup_{\t\in\mathcal{S}_{n}}Z^\alpha(\t)>u} \sim\sum_{j\in \mathcal{N}}\pk{\sup_{\t\in\Pi_{\delta, j}}Z^\alpha(\t)>u}\sim \mathfrak{m}\Psi(u), \ u\rw\IF.
\EQNY
This completes the proof.
\QED
\section{Proof of Proposition \ref{ex02} }\label{chiproof}
Observe that for $0<\epsilon<\pi/4$
\begin{align}
\mathbb{P}\left(\sup_{(\vk{\theta}, t)\in E_{1,\epsilon}}Z(\vk{\theta},t)>u\right)\leq \mathbb{P}\left(\sup_{(\vk{\theta}, t)\in E}Z(\vk{\theta},t)>u\right)\leq \sum_{i=1}^{3}\mathbb{P}\left(\sup_{(\vk{\theta}, t)\in E_{i,\epsilon}}Z(\vk{\theta},t)>u\right),
\end{align}
where
$$E_{1,\epsilon}=[\epsilon,\pi-\epsilon]^{n-2}\times [0,2\pi-\epsilon)\times[0,\epsilon],~ E_{2,\epsilon}=[0,\pi]^{n-2}\times [0,2\pi)\times[\epsilon, 1],~E_{3,\epsilon}=E/(E_{1,
\epsilon}\cup E_{2,\epsilon}).
$$
We will first apply Theorem \ref{MThm2} to obtain the asymptotics over $E_{1,\epsilon}$ and then show that the asymptotics over $E_{2,\epsilon}$ and $E_{3,\epsilon}$ are negligible as $u\to\infty $ and $\epsilon\to 0$.

{\underline{The asymptotics over $E_{1,\epsilon}$}}.  To this end, we next analyze the variance and correlation of $Z$.
By (\ref{X}), we have
\begin{align}\label{varZ}\sigma_Z(\vk{\theta},t)=\frac{1}{1+bt^{\alpha}},~t\in [0,1].
\end{align}
Hence $\sigma_Z(\vk{\theta},t)$ attains its maximum equal to 1 at $[0,\pi]^{n-2}\times[0,2\pi)\times\{0\}$ and
$$\lim_{\delta\downarrow 0}\sup_{\vk{\theta}\in [0,\pi]^{n-2}\times[0,2\pi), 0<t<\delta }\left|\frac{1-\sigma_Z(\vk{\theta},t)}{bt^{\alpha}}-1\right|=1.$$
This implies that {\bf A2} is satisfied. For {\bf A1},
 by (\ref{Y}), we have
\begin{align*}1-Corr(Z(\vk{\theta},t), Z(\vk{\theta}',t'))&\sim aVar(Y(t)-Y(t'))+\frac{1}{2}\sum_{i=1}^{n}(v_i(\vk{\theta})-v_i(\vk{\theta}'))^2\\
& \sim aVar(Y(t)-Y(t'))+\frac{(\theta_1-\theta_1')^2}{2}+\frac{1}{2}\sum_{i=2}^{n-1}\left(\prod_{j=1}^{i-1}\sin(\theta_j)\right)^2  (\theta_i-\theta_i')^2,
\end{align*}
as  $(\vk{\theta},t), (\vk{\theta}',t')\in E$ and $|t-t'|, |\vk{\theta}-\vk{\theta}'|\to 0$.
Let \begin{align}\label{W3}
W(\vk{\theta}, t)=\sum_{i=1}^{n-1}B_i^2(\theta_i)+\sqrt{a}Y(t),~\vk{\theta}\in \mathbb{R}^{n-1}\times \mathbb{R}^+,
\end{align} where $B_i^2$ are independent fractional Brownian motions with index $2$ and $Y$ is the self-similar Gaussian process given in (\ref{Y}) that is independent of $B_i^2$. Moreover, we denote $\vk{a}(\vk{\varphi})=(a_1(\vk{\varphi}),
\dots, a_{n-1}(\vk{\varphi})),~\vk{\varphi}\in[0,\pi]^{n-2}\times [0,2\pi)$ with  $$a_1(\vk{\varphi})=\frac{1}{\sqrt{2}}~ \text{and}~ a_i(\vk{\varphi})=\frac{1}{
\sqrt{2}}\prod_{j=1}^{i-1}\sin(\varphi_j),~i=2,\dots,n-1.$$
It follows that for $0<\epsilon<\pi/4$
$$\lim_{\delta\downarrow 0}\sup_{\vk{\varphi}\in [\epsilon,\pi-\epsilon]^{n-2}\times [0,2\pi)}\sup_{(\vk{\theta},t), (\vk{\theta}',t')\in E, |(\vk{\theta},t)-(\vk{\varphi},0)|, |(\vk{\theta}',t')-(\vk{\varphi},0)|<\delta }\left|\frac{1-Corr(Z(\vk{\theta},t), Z(\vk{\theta}',t'))}{\E{\left(W(\vk{a}(\vk{\varphi})\vk{\theta},t)-W(\vk{a}(\vk{\varphi})\vk{\theta}',t')\right)^2}}-1\right|=0.$$
By the fact that
\begin{align}\label{WY}
Var(W(\vk{\theta},t)-W(\vk{\theta}',t'))=aVar(Y(t)-Y(t'))+\sum_{i=1}^{n-1}(\theta_i-\theta_i')^2,
\end{align}
we know that $W(\vk{\theta},t)$ is homogeneous with respect to $\vk{\theta}$ if $t$ is fixed.  This implies that (\ref{Rr}) holds with $W$ defined in (\ref{W3}).

Moreover, by self-similarity of $Y$ and (\ref{WY}) we have
\begin{align*}
Var(W(u^{-1}\vk{\theta},u^{-2/\alpha}t)-W(u^{-1}\vk{\theta}', u^{-2/\alpha}t'))=u^{-2}Var(W(\vk{\theta},t)-W(\vk{\theta}',t')),
\end{align*}
showing that (\ref{W}) holds with $\alpha_i=2, ~i=1,\dots, n-1$ and $\alpha_n=\alpha$.
In addition, by {\bf B1-B2}, there exists $d>0$ such that for $|\vk{\theta},t)-(\vk{\theta}',t')|<\delta$ with $(\vk{\theta},t), (\vk{\theta}',t')\in E_{1,\epsilon}$,
$$\mathbb{Q}_1\sum_{i=1}^{n-1}(\theta_i-\theta_i')^2\leq 1-Corr(Z(\vk{\theta},t)\leq \mathbb{Q}_2\left(|t-t'|^{\alpha}+ \sum_{i=1}^{n-1}(\theta_i-\theta_i')^2\right).$$ Hence (\ref{boundrr1}) is confirmed. Moreover, (\ref{reql}) is clearly satisfied over $E_{1,\epsilon}$. Therefore, {\bf A1} is verified for $Z$ over $E_{1,\epsilon}$.
Consequently, it follows from Theorem \ref{MThm2} that, as $u\to\infty$,
$$\mathbb{P}\left(\sup_{(\vk{\theta}, t)\in E_{1,\epsilon}}Z(\vk{\theta},t)>u\right)\sim \mathcal{H}_W^{bt^{\alpha}}\int_{\vk{\theta}\in [\epsilon,
\pi-\epsilon]^{n-2}\times [0,2\pi-\epsilon]} 2^{-(n-1)/2}\prod_{i=1}^{n-1}|\sin(\theta_i)|^{n-i-1}d \theta_1\dots d\theta_{n-1}u^{n-1}\Psi(u),$$
where $W$ is given in (\ref{W3}).

{\underline{Upper bound for the asymptotics over $E_{2,\epsilon}$}}. By (\ref{varZ}), there exists $0<\delta<1$ such that
$$\sup_{(\vk{\theta},t)\in E_{2,\epsilon}}Var(Z(\vk{\theta},t))\leq 1-\delta.$$
It follows from Borell-TIS inequality that
\begin{align*}
\mathbb{P}\left(\sup_{(\vk{\theta}, t)\in E_{2,\epsilon}}Z(\vk{\theta},t)>u\right)\leq \exp\left(-\frac{\left(u-\E{\sup_{(\vk{\theta}, t)\in E_{2,\epsilon}}Z(\vk{\theta},t)}\right)^2}{2(1-\delta)}\right)=o(u^{n-1}\Psi(u)),~u\to\infty.
\end{align*}
{\underline{Upper bound for the asymptotics over $E_{3,\epsilon}$}}. Direct calculation shows that
$$1-Corr(Z(\vk{\theta},t)\leq \mathbb{Q}_2\left(|t-t'|^{\alpha}+ \sum_{i=1}^{n-1}(\theta_i-\theta_i')^2\right)$$
holds for $(\vk{\theta},t), (\vk{\theta}',t')\in E_{3,\epsilon}$. Define $U_3(\vk{\theta},t), (\vk{\theta},t)\in \R^n$ to be a centered homogeneous Gaussian field with continuous trajectories, unit variance and correlation function $r_{U_3}(\vk{\theta},t, \vk{\theta}',t')$ satisfying
\BQNY
r_{U_3}(\vk{\theta},t, \vk{\theta}',t')=1-\exp\left(-2\mathbb{Q}_2\left(|t-t'|^{\alpha}+ \sum_{i=1}^{n-1}(\theta_i-\theta_i')^2\right)\right).
\EQNY
By Slepian's inequality and Theorem \ref{MThm2}, we have
\begin{align*}
\mathbb{P}\left(\sup_{(\vk{\theta}, t)\in E_{3,\epsilon}}Z(\vk{\theta},t)>u\right)\leq \mathbb{P}\left(\sup_{(\vk{\theta}, t)\in E_{3,\epsilon}}\frac{U_3(\vk{\theta},t)}{1+bt^{\alpha}}>u\right)\leq \mathbb{Q}v_n(E_{3,\epsilon})u^{n-1}\Psi(u),~u\to\infty.
\end{align*}
Noting that $\lim_{\epsilon\to 0}v_n(E_{3,\epsilon})=0$, combination of the above asymptotics and upper bounds leads to
$$\mathbb{P}\left(\sup_{(\vk{\theta}, t)\in E}Z(\vk{\theta},t)>u\right)\sim \mathcal{H}_W^{bt^{\alpha}}\int_{\vk{\theta}\in [0,
\pi]^{n-2}\times [0,2\pi)} 2^{-(n-1)/2}\prod_{i=1}^{n-1}|\sin(\theta_i)|^{n-i-1}d \theta_1\dots d\theta_{n-1}u^{n-1}\Psi(u),~u\to\infty.$$
By the fact that
$$\int_{\vk{\theta}\in [0,
\pi]^{n-2}\times [0,2\pi)} \prod_{i=1}^{n-1}|\sin(\theta_i)|^{n-i-1}d \theta_1\dots d\theta_{n-1}=\frac{2\pi^{n/2}}{\Gamma(n/2)},$$
and
$\mathcal{H}_W^{bt^{\alpha}}=\mathcal{P}_{\sqrt{a}Y}^{b}(\mathcal{H}_{B^2})^{n-1}
=\mathcal{P}_{Y}^{a^{-1}b}\pi^{-(n-1)/2},$ where we used that $\mathcal{H}_{B^2}=\pi^{-1/2},$
we have
$$\mathbb{P}\left(\sup_{(\vk{\theta}, t)\in E}Z(\vk{\theta},t)>u\right)\sim
\frac{2^{\frac{3-n}{2}}\sqrt{\pi}}{\Gamma(n/2)}\mathcal{P}_{Y}^{a^{-1}b}u^{n-1}\Psi(u),~u\to\infty.$$
  \QED
\section{Appendix: proof of Remark \ref{prop41}}
i) For the case $1\leq \mathfrak{m}\leq n$,
we first show that $\mathcal{H}_W\geq 1$.
Recall that $\mathcal{N}_0=\{i\in\mathcal{N}, i<k^*\}$, $\MNB=\{i: a_i<1, i=1\ldot n+1\}$ and $\widetilde{\x}=(x_i)_{i\in \mathcal{N}_0 }$.

 For $x_i=0, ~i\in \MNB$, by the definition of $W$ in (\ref{defineW}), we have  $$\left\{\sqrt{2}W(\x)-\sum_{i=1, i\neq k^*}^{n+1}x_i,~~ \widetilde{\x}\in[0,\lambda]^{\mathfrak{m}-1}\right\}\laweq \left\{\sum_{i\in \mathcal{N}_0}\sqrt{2}B_i(x_i)
-\sum_{i\in\mathcal{N}_0}x_i,~~ \widetilde{\x}\in[0,\lambda]^{\mathfrak{m}-1}\right\}.$$
Hence
\BQNY
\mathcal{H}_W\geq
\lim_{\lambda\rw\IF}\frac{1}{\lambda^{\mathfrak{m}-1}}
\E{\sup_{\widetilde{\x}\in[0,\lambda]^{\mathfrak{m}-1}}e^{\sum_{i\in \mathcal{N}_0}\sqrt{2}B_i(x_i)
-\sum_{i\in\mathcal{N}_0}x_i}}
=\prod_{i\in \mathcal{N}_0}\mathcal{H}_{B_i},
\EQNY
where $H_{B_i}$ is defined in (\ref{H1}). Note that $\mathcal{H}_{B_i}=1$, see e.g., \cite{Pit96} (or \cite{GP2018}). Therefore,
$\mathcal{H}_W\geq 1.$
We next derive the upper bound of $\mathcal{H}_W$ for $1\leq \mathfrak{m}\leq n$, for which we need to use the notation  in the proof of ii) of Proposition \ref{ThmM1}, e.g., $Y$ and $\widetilde{\mathcal{S}}_n(\delta)$.
For $\delta\in (0, \frac{1}{(n+1)^2})$, let
$$A(\delta)=\{\x: \widetilde{\x}\in B(\delta),~0\leq x_i\leq \frac{\delta}{(n+1)^2}, \text{if} \ i\in\MNB\},$$
where $B(\delta)=\prod_{i=1}^{m-1}[2i\delta, (2i+1)\delta]$. Clearly, $A(\delta)\subset\widetilde{\mathcal{S}}_n(\delta)$.  Moreover, by (\ref{yrr}) we have that for any $\epsilon>0$, there exists $\delta \in (0, \frac{1}{(n+1)^2})$ such that  for any $\x,\y\in A(\delta)$,
\BQNY
1-r_Y(\x,\y)\leq (n+\epsilon)\underset{i\neq k^*}{\sum_{i=1}^{n+1}}\abs{x_i-y_i}.
\EQNY
Let us introduce a centered homogeneous Gaussian fields  $U_4(\x)$, $\x\in[0,\infty)^{n}$ with continuous trajectories,  unit variance and correlation functions
$$
r_{U_4}(\x,\y)=\exp\LT(-\E{\LT(W_4(\x)-W_4(\y)\RT)^2}\RT),~\text{with}~
W_4(\x)=\sqrt{n+\epsilon}\underset{i\neq k^*}{\sum_{i=1}^{n+1}}  B_i(x_i),$$
where $B_i, i=1,\dots, k^*-1, k^*+1, n+1$ are iid standard Brownian motions.
By (\ref{eqcona1}) and Slepian's inequality, we have, for $0<\epsilon<1$,
\BQNY
\pk{\sup_{\x\in A(\delta)}\frac{U_4(\x)}{1+\sum_{i\in\MNB}\frac{1-a_i^2}{2-\epsilon}x_i}>u}\geq \pk{\sup_{\x\in A(\delta)} Y(\x)>u}.
\EQNY
Analogously to (\ref{asyp1}), we have
\BQNY
\pk{\sup_{\x\in A(\delta)}Y(\x)>u}\sim v_{\mathfrak{m}-1}\LT(B(\delta)\RT)\mathcal{H}_Wu^{2(\mathfrak{m}-1)}\Psi(u),
\EQNY
and
\BQNY
\pk{\sup_{\x\in A(\delta)}\frac{U_4(\x)}{1+\sum_{i\in\MNB}\frac{1-a_i^2}{2+\epsilon}x_i}>u}\sim v_{\mathfrak{m}-1}\LT(B(\delta)\RT)\mathcal{H}_{W_4} u^{2(\mathfrak{m}-1)}\Psi(u),
\EQNY
where
\begin{align*}
  \mathcal{H}_{W_4}&=\lim_{\lambda\rw\IF}\frac{1}{\lambda^{\mathfrak{m}-1}}
\E{\sup_{\x\in[0,\lambda]^n}e^{ \sqrt{2(n+\epsilon)}\underset{i\neq k^*}{\sum_{i=1}^{n+1}}B_i(x_i)
-(n+\epsilon)\underset{i\neq k^*}{\sum_{i=1}^{n+1}}x_i-\sum_{i\in\MNB}\frac{1-a_i^2}{2-\epsilon}x_i}}\\
&=(n+\epsilon)^{\mathfrak{m}-1}\left(\prod_{i
\in \mathcal{N}_0}\mathcal{H}_{B_i}\right)\prod_{i\in \mathcal{N}^c} \mathcal{P}_{B_i}^{\frac{1-a_i^2}{(2-\epsilon)(n+\epsilon)}},
\end{align*}
with $\mathcal{P}_{B_i}^c$ for $c>0$ being defined in (\ref{H1}).
Using the fact that $\mathcal{H}_{B_i}=1$  and for $c>0$ (see, e.g. \cite{GP2018}),
$
\mathcal{P}_{B_i}^c=1+\frac{1}{c},
$
we have
$$\mathcal{H}_{W_4}=(n+\epsilon)^{\mathfrak{m}-1}\prod_{i\in\MNB}\left(1+\frac{(2+\epsilon)(n+\epsilon)}{1-a_i^2}\right).$$
Hence
$$\mathcal{H}_{W}\leq \mathcal{H}_{W_4}=(n+\epsilon)^{\mathfrak{m}-1}\prod_{i\in\MNB}\left(1+\frac{(2+\epsilon)(n+\epsilon)}{1-a_i^2}\right).$$
We establish the claim by letting $\epsilon\to 0$.

ii) If $\mathfrak{m}=n+1$, we have $\mathcal{N}_0=\{1,\dots, n\}$ and
\BQNY
\mathcal{H}_W=
\lim_{\lambda\rw\IF}\frac{1}{\lambda^{n}}
\E{\sup_{\widetilde{\x}\in[0,\lambda]^{n}}e^{\sum_{i\in \mathcal{N}_0}\sqrt{2}B_i(x_i)
-\sum_{i\in\mathcal{N}_0}x_i}}
=\prod_{i\in \mathcal{N}_0}\mathcal{H}_{B_i}=1.
\EQNY
\QED

\section*{Acknowledgments}
We thank Enkelejd Hashorva for many stimulating comments that highly
improved the content of this contribution. We thank Lanpeng Ji for some useful discussions.
Support from SNSF Grant 200021-175752/1 is kindly acknowledged.  KD
was partially supported by NCN Grant No 2018/31/B/ST1/00370
(2019-2022).
Long Bai is supported by  National Natural Science Foundation of China Grant no. 11901469 and The Natural Science Foundation
of the Jiangsu Higher Education Institutions of China grant no. 19KJB110022.

\bibliographystyle{ieeetr}

\COM{\newcommand{\nosort}[1]{}\def\polhk#1{\setbox0=\hbox{#1}{\ooalign{\hidewidth
  \lower1.5ex\hbox{`}\hidewidth\crcr\unhbox0}}}
  \def\polhk#1{\setbox0=\hbox{#1}{\ooalign{\hidewidth
  \lower1.5ex\hbox{`}\hidewidth\crcr\unhbox0}}} \def\cprime{$'$}
  \def\cprime{$'$} \def\cprime{$'$}
  \def\lfhook#1{\setbox0=\hbox{#1}{\ooalign{\hidewidth
  \lower1.5ex\hbox{'}\hidewidth\crcr\unhbox0}}}
  \def\polhk#1{\setbox0=\hbox{#1}{\ooalign{\hidewidth
  \lower1.5ex\hbox{`}\hidewidth\crcr\unhbox0}}} \def\cprime{$'$}
}

\bibliography{GUE}

\newcommand{\nosort}[1]{}\def\polhk#1{\setbox0=\hbox{#1}{\ooalign{\hidewidth
  \lower1.5ex\hbox{`}\hidewidth\crcr\unhbox0}}}
  \def\polhk#1{\setbox0=\hbox{#1}{\ooalign{\hidewidth
  \lower1.5ex\hbox{`}\hidewidth\crcr\unhbox0}}} \def\cprime{$'$}
  \def\cprime{$'$} \def\cprime{$'$}
\begin{thebibliography}{10}

\bibitem{Pit96}
V.~I. Piterbarg, {\em Asymptotic methods in the theory of {G}aussian processes
  and fields}, vol.~148 of {\em Translations of Mathematical Monographs}.
\newblock American Mathematical Society, 1996.

\bibitem{AdlerTaylor}
R.~Adler and J.~Taylor, {\em Random fields and geometry}.
\newblock Springer Monographs in Mathematics, New York: Springer, 2007.

\bibitem{Lif13}
M.~A. Lifshits, {\em {G}aussian random functions}, vol.~322.
\newblock Springer Science \& Business Media, 2013.

\bibitem{PiP81}
V.~I. Piterbarg and V.~Prisyazhnyuk, ``The exact asymptotics for the
  probability of large span of a {G}aussian stationary process,'' {\em Teoriya
  Veroyatnostei i ee Primeneniya}, vol.~26, no.~3, pp.~480--495, 1981.

\bibitem{Liu16}
P.~Liu, ``Extremes of {G}aussian random fields with maximum variance attained
  over smooth curves,'' {\em arXiv preprint arXiv:1612.07780}, 2016.

\bibitem{Liu18}
D.~Cheng and P.~Liu, ``Extremes of spherical fractional {B}rownian motion,''
  {\em Extremes}, vol.~22, no.~3, pp.~433--457, 2019.

\bibitem{DHL16}
K.~D{\polhk{e}}bicki, E.~Hashorva, and L.~Ji, ``{Extremes of a class of
  nonhomogeneous {G}aussian random fields},'' {\em Ann. Probab.}, vol.~44,
  no.~2, pp.~984 -- 1012, 2016.

\bibitem{Chan2006}
H.~P. Chan and T.~L. Lai, ``Maxima of asymptotically {G}aussian random fields
  and moderate deviation approximations to boundary crossing probabilities of
  sums of random variables with multidimensional indices,'' {\em Ann. Probab.},
  vol.~34, no.~1, pp.~80--121, 2006.

\bibitem{Adler1986}
R.~Adler and L.~Brown, ``Tail behavior for suprema of empirical processes,''
  {\em Ann. Probab.}, vol.~14, no.~1, pp.~1--30, 1986.

\bibitem{Pit21}
V.~I. Piterbarg, ``High excursion probabilities for {G}aussian fields on smooth
  manifolds,'' {\em arXiv preprint arXiv:2108.07473}, 2021.

\bibitem{Bar01}
Y.~Baryshnikov, ``Gues and queues,'' {\em Probab. Theory Relat. Fields},
  vol.~119, no.~2, pp.~256--274, 2001.

\bibitem{Oco02}
N.~O'Connell, ``Random matrices, non-colliding processes and queues,'' {\em
  S{\'e}minaire de probabilit{\'e}s de Strasbourg}, vol.~36, pp.~165--182,
  2002.

\bibitem{GlW91}
P.~W. Glynn and W.~Whitt, ``Departures from many queues in series,'' {\em Ann.
  Appl. Probab.}, pp.~546--572, 1991.

\bibitem{GTW01}
J.~Gravner, C.~A. Tracy, and H.~Widom, ``Limit theorems for height fluctuations
  in a class of discrete space and time growth models,'' {\em J. Stat. Phys.},
  vol.~102, no.~5, pp.~1085--1132, 2001.

\bibitem{Lindgren1980a}
G.~Lindgren, ``Extreme values and crossing for the chi-square processes and
  other functions of multidimensional {G}aussian process, with reliability
  applications,'' {\em Adv. Appl. Probab.}, vol.~12, pp.~746--774, 1980.

\bibitem{Pitchi1994}
V.~I. Piterbarg, ``High excursions for nonstationary generalized chi-square
  processes,'' {\em Stoch. Process. Appl.}, vol.~53, pp.~307--337, 1994.

\bibitem{EnkelejdJi2014Chi}
E.~Hashorva and L.~Ji, ``Piterbarg theorems for chi-processes with trend,''
  {\em Extremes}, vol.~18, pp.~37--64, 2015.

\bibitem{PL2015}
P.~Liu and L.~Ji, ``Extremes of chi-square processes with trend,'' {\em Probab.
  Math. Statist.}, vol.~36, pp.~1--20, 2016.

\bibitem{LJ2017}
P.~Liu and L.~Ji, ``Extremes of locally stationary chi-square processes with
  trend,'' {\em Stoch. Process. Appl.}, vol.~127, pp.~497--525, 2017.

\bibitem{Bai2021}
L.~Bai and D.~Kalaj, ``Approximation of {K}olmogorov-{S}mirnov test
  statistics,'' {\em Stochastics}, vol.~93, pp.~993--1027, 2021.

\bibitem{Sri93}
R.~Srinivasan, ``Queues in series via interacting particle systems,'' {\em
  Math. Oper. Res.}, vol.~18, no.~1, pp.~39--50, 1993.

\bibitem{Houdre2003}
C.~Houdr\'e and J.~Villa, ``An example of infinite dimensional quasi-helix,''
  {\em in: Stochastic Models (Mexico City, 2002), Contemp. Math. 336, Amer.
  Math. Soc.}, pp.~195--202, 2003.

\bibitem{Lei2009}
P.~Lei and D.~Nualart, ``A decomposition of the bifractional {B}rownian motion
  and some applications,'' {\em Statist. Probab. Lett.}, vol.~79, pp.~619--624,
  2009.

\bibitem{Bojdecki2004}
T.~Bojdecki and A.~Gorostiza, L. G.and~Talarczyk, ``Sub-fractional {B}rownian
  motion and its relation to occupation times,'' {\em Statist. Probab. Lett.},
  vol.~69, pp.~405--419, 2004.

\bibitem{Dzhaparidze2004}
K.~Dzhaparidze and H.~Van~Zanten, ``A series expansion of fractional {B}rownian
  motion,'' {\em Probab. Theory Relat. Fields}, vol.~130, pp.~39--55, 2004.

\bibitem{Li2004}
W.~Li and Q.~Shao, ``Lower tail probabilities for {G}aussian processes,'' {\em
  Ann. Probab.}, vol.~32, pp.~216--242, 2004.

\bibitem{Debic2020}
K.~D{\polhk{e}}bicki and K.~Tabi\'s, ``Pickands-{P}iterbarg constants for
  self-similar {G}aussian processes,'' {\em Probab. Math. Statist.}, vol.~40,
  no.~2, pp.~297--315, 2020.

\bibitem{Fatalov1993}
V.~Fatalov, ``Asymptotics of large deviation probabilities for {G}aussian
  fields: Applications.,'' {\em Izvestiya Natsionalnoi Akademii Nauk Armenii},
  vol.~28, pp.~25--51, 1993.

\bibitem{Uniform2016}
K.~D{\polhk{e}}bicki, E.~Hashorva, and P.~Liu, ``Uniform tail approximation of
  homogenous functionals of {G}aussian fields,'' {\em Adv. Appl. Probab.},
  vol.~49, no.~4, pp.~1037--1066, 2017.

\bibitem{GP2018}
L.~Bai, K.~D{\polhk{e}}bicki, E.~Hashorva, and L.~Luo, ``On generalised
  {P}iterbarg constants,'' {\em Methodol. Comput. Appl. Probab.}, 2017.

\end{thebibliography}


\begin{thebibliography}{1}

\bibitem{AdlerTaylor}
R.~Adler and J.~Taylor, {\em Random fields and geometry}.
\newblock Springer Monographs in Mathematics, New York: Springer, 2007.
\bibitem{Pickands1969}
III. J. Pickands, ``Upcrossing probabilities for stationary {G}aussian processes," {\em Trans. Amer. Math. Soc. }, vol.~145, pp.~ 51--73, 1969.

\bibitem{Pickands1969b}
III. J. Pickands, ``Asymptotic properties of the maximum in a stationary Gaussian process," {\em Trans. Amer. Math. Soc. }, vol.~145, pp.~ 75--86, 1969.

\bibitem{Piterbarg1978}
V. I. Piterbarg and V. P. Prisja\v{z}njuk, ``Asymptotic behavior of the probability of a large excursion for a
nonstationary {G}aussian process'' {\em Teor. Verojatnost. i Mat. Statist.}, vol.~18, pp.~121--134, 1978.

\bibitem{lp2018}
L.~Bai, ``Extremes of Lp-norm of vector-valued {G}aussian
processes with trend,'' {\em Stochastics:
An International Journal of Probability and Stochastic Processes}, vol.~90, pp.~1111--1144, 2018.

\bibitem{GP2018}
L.~Bai, K.~D{\polhk{e}}bicki, E.~Hashorva, and L.~Luo, ``On generalised {P}iterbarg constants,'' {\em Methodol Comput Appl Probab}, vol.~20, pp.~137--164, 2018.


\bibitem{Uniform2016}
K.~D{\polhk{e}}bicki, E.~Hashorva, and P.~Liu, ``Uniform tail approximation of
  homogenous functionals of {G}aussian fields,'' {\em Advances in Applied
  Probability}, vol.~49, pp.~1037--1066, 2017.

\bibitem{Pit96}
V.~Piterbarg, {\em Asymptotic methods in the theory of {G}aussian processes and
  fields}, vol.~148 of {\em Translations of Mathematical Monographs}.
\newblock Providence, RI: American Mathematical Society, 1996.



\bibitem{GennaBorell}
G.~Samorodnitsky, ``Probability tails of {G}aussian extrema,'' {\em Stochastic
  Process. Appl.}, vol.~38, no.~1, pp.~55--84, 1991.

\end{thebibliography}

\end{document}